\documentclass[11pt]{amsart}
\usepackage{latexsym}
\usepackage{amsfonts,amssymb,amsmath,amsthm,amscd}
\usepackage[mathscr]{eucal}

\newcommand{\af}{\mbox{${\mathfrak a}$}}

\newcommand{\g}{\mbox{${\mathfrak g}$}}
\newcommand{\h}{\mbox{${\mathfrak h}$}}
\newcommand{\kf}{\mbox{${\mathfrak k}$}}

\newcommand{\m}{\mbox{${\mathfrak m}$}}

\newcommand{\p}{\mbox{${\mathfrak p}$}}

\newcommand{\C}{\mbox{${\mathbb C}$}}

\newcommand{\PP}{\mbox{${\mathbb P}$}}

\newcommand{\R}{\mbox{${\mathbb R}$}}
\newcommand{\Z}{\mbox{${\mathbb Z}$}}

\newcommand{\tr}{{\rm tr}}
\newcommand{\ric}{{\rm Ric}}

\makeatletter
\def\numberwithin#1#2{\@ifundefined{c@#1}{\@nocnterrr}{%
  \@ifundefined{c@#2}{\@nocnterr}{%
  \@addtoreset{#1}{#2}%
  \toks@\expandafter\expandafter\expandafter{\csname the#1\endcsname}%
  \expandafter\xdef\csname the#1\endcsname
    {\expandafter\noexpand\csname the#2\endcsname
     .\the\toks@}}}}
\makeatother
\numberwithin{equation}{section}

\newtheorem{thm}[equation]{Theorem}
\newtheorem{lemma}[equation]{Lemma}
\newtheorem{prop}[equation]{Proposition}

\newtheorem{cor}[equation]{Corollary}
\newtheorem{ex}[equation]{Example}
\newenvironment{example}{\begin{ex} \em}{\end{ex}}
\newtheorem{rem}[equation]{Remark}
\newenvironment{rmk}{\begin{rem} \em}{\end{rem}}

\begin{document}

\title{On Ricci Solitons of Cohomogeneity One}
\author{Andrew S. Dancer}
\address{Jesus College, Oxford University, OX1 3DW, United Kingdom}
\email{dancer@maths.ox.ac.uk}
\author{McKenzie Y. Wang}
\address{Department of Mathematics and Statistics, McMaster
University, Hamilton, Ontario, L8S 4K1, Canada}
\email{wang@mcmaster.ca}
\thanks{partially supported by NSERC Grant No. OPG0009421}

\date{revised \today}

\begin{abstract}
We analyse some properties of  the cohomogeneity one Ricci soliton
equations, and use Ans{\" a}tze of cohomogeneity one type to produce new
explicit examples of complete K\"ahler Ricci solitons of
expanding, steady and shrinking types. These solitons are foliated by
hypersurfaces which are circle bundles over a product of Fano
K\"ahler-Einstein manifolds or over coadjoint orbits of a compact
connected semisimple Lie group.
\end{abstract}

\maketitle


\bigskip
\setcounter{section}{-1}

\section{\bf Introduction}

A {\em Ricci soliton} is a solution $(g,X)$, where $g$ is a complete Riemannian
metric and $X$ is a complete vector field on a manifold $M$, to the equation:
\begin{equation} \label{LieRS}
{\ric}(g) + \frac{1}{2} \,{\mathscr L}_{X} g + \frac{\epsilon}{2} \, g = 0,
\end{equation}
where $\epsilon$ is a real constant.
The significance of such a solution is that it generates a family
of metrics that evolves in a particularly simple way under the Ricci flow
\begin{equation} \label{ricciflow}
\frac{\partial g_{\tau}}{\partial \tau} = - 2 {\ric}(g_{\tau}).
\end{equation}
For we may now define a 1-parameter family of vector fields
\[
Y_{\tau} = \frac{1}{1 + \epsilon \tau} \, X
\]
and integrate these to a $1$-parameter family of diffeomorphisms $\psi_{\tau}$
on $M$. Then $g_{\tau} = (1 + \epsilon \tau) \psi_{\tau}^* g$ defines
a solution to the Ricci flow which evolves just by diffeomorphisms
and homotheties. The soliton is called {\em steady}, {\em expanding,} or
{\em shrinking} depending on whether $\epsilon$ is zero, positive, or negative.
Steady solitons therefore give examples of {\em eternal} Ricci flows, i.e., ones
defined for all $\tau \in (-\infty, \infty)$. Expanding and shrinking
solitons give respectively {\em immortal} solutions on $(-\frac{1}{\epsilon},+\infty)$
 and {\em ancient} solutions on $(-\infty,\frac{1}{|\epsilon|}).$

Note that equation (\ref{LieRS}) may be written instead as
\begin{equation} \label{formRS}
{\ric}(g) + \delta^* \omega  + \frac{\epsilon}{2} \, g =0
\end{equation}
where $\omega$ is the 1-form dual to $X$ via the metric $g$ and $\delta^*$
is the symmetrized covariant derivative.
A particularly important class of solutions is obtained if we take
$\omega$ to be exact, or equivalently take $X$ to be the gradient of
a smooth function $u$. In this case the Ricci soliton equation becomes
\begin{equation} \label{HessRS}
{\ric}(g) + {\rm Hess}(u) + \frac{\epsilon}{2}\, g =0,
\end{equation}
where $\rm Hess$ denotes the Hessian.  A solution $(g, u)$ of this equation
is called a {\em Ricci soliton of gradient type}.

Ricci solitons are of course generalisations of Einstein metrics, for
if $X$ is the zero vector field we recover the Einstein equation from
(\ref{LieRS}). It is natural therefore to ask whether techniques for
producing Einstein metrics can be adapted to produce examples of Ricci
solitons. One such approach is to look for solutions with large
symmetry group. Recall that Perelman's no breathers theorems (cf \S 2 and 3
of \cite{Per}) imply that on a compact manifold all Ricci solitons are
of gradient type. Several authors, e.g., \cite{ELM} have then observed that
this fact immediately implies that compact Ricci solitons with constant
scalar curvature must be {\em trivial}, i.e., Einstein.
For if we take the trace of (\ref{HessRS})
we obtain
\[
R + \Delta u + \frac{n \epsilon}{2} =0,
\]
where $R$ is the scalar curvature and $n$ is the dimension of $M$.
If $R$ is constant,
 we can integrate over the manifold to
deduce that $R + \frac{n \epsilon}{2} =0,$ so $u$ is constant and the soliton
is trivial.
In particular, compact  homogeneous solitons must be trivial.
By contrast, nontrivial {\em noncompact} homogeneous Ricci
solitons exist, and there is a beautiful correspondence between
left-invariant Einstein metrics on solvable Lie groups and
homogeneous Ricci solitons on their nilradicals (cf \cite{La}, \cite{FDC}).

The next step up in complexity from homogeneous metrics are those of
{\em cohomogeneity one}, that is, those where a group acts isometrically
with generic orbit of codimension one. Curvature equations are therefore
reduced to ordinary differential equations in a variable transverse to the
orbits. The foundational work here is due to B\'erard Bergery \cite{BB},
who developed the mathematical framework of Page's metric \cite{Pa}
and produced further new examples of Einstein metrics of cohomogeneity one.
We remark here that the same equations may arise in any situation where
the metric on the hypersurfaces depends on a single transverse variable,
whether or not the hypersurfaces are homogeneous (see Remark \ref{gen}). Indeed many of
B\'erard Bergery's examples were already of this form.

Recall that a Ricci soliton is said to be {\em K\"ahler} if in addition
there is a complex structure for which $g$ is K\"ahler and $X$ is an
infinitesimal automorphism. By the work of Tian and Zhu \cite{TZ1}, \cite{TZ2},
K\"ahler-Ricci solitons on {\em compact} complex manifolds
are unique up to holomorphic transformations. Several authors
\cite{Ko}, \cite{Ca1}, \cite{Ca2}, \cite{ChV}, \cite{G}, \cite{PTV},
\cite{FIK}, \cite{Yan} have produced cohomogeneity one type K\"ahler-Ricci
solitons where the hypersurfaces are circle bundles over a Fano
K\"ahler-Einstein space. The hypersurfaces are equipped with metrics
such that the bundle projection becomes a Riemannian submersion with
totally geodesic fibres. We also note that X.-J. Wang and Xiaohua Zhu
have shown that compact toric Fano manifolds always admit K\"ahler Ricci
solitons \cite{WZ}, and their result has been generalised in \cite{PS2} to
Fano bundles over a generalised flag manifold with compact toric Fano fibres.

In this paper we formulate a general approach to Ricci solitons of
cohomogeneity one and use it, in the K\"ahler setting, to unify and
generalise all the above cohomogeneity one type examples. We
particularly focus on the situation where the hypersurfaces are
generic circle fibrations over arbitrary compact homogeneous K\"ahler
manifolds (see \S 4), or certain circle bundles over an arbitrary
product of Fano K\"ahler-Einstein manifolds (see \S 3). In the latter
case, the Fano manifolds need not have any isometries, so
the resulting K\"ahler manifolds need not have more than a circle of
isometries. We consider both the compact and non-compact cases, and
for each hypersurface we analyse the different blow-downs which form a
smooth compact end. We note that the corresponding Fano
K\"ahler-Einstein case was considered by \cite{Sa}, \cite{KS1},
\cite{KS2}, \cite{PS1}, while the complete non-compact
K\"ahler-Einstein case was examined in \cite{DW} and \cite{WW}, (cf
\cite{Wa}, Thms 3.1 and 3.2 for more general blow-downs than those in
\cite{WW}).

In more detail, the layout of the paper is as follows. In \S 1 we develop
the general formalism for cohomogeneity one Ricci solitons, and write down the
resulting system of ordinary differential equations. \S 2 is devoted
to proving an analogue for Ricci solitons of
a theorem of A. Back  concerning the cohomogeneity one Einstein equations
\cite{Ba}. Namely, we show that provided a special orbit of dimension
strictly smaller than that of the principal orbits is present and a sufficient
amount of smoothness of the metric and $1$-form is established,
the full Ricci soliton equations on a cohomogeneity one manifold actually
follow from a smaller set of equations which includes the
components of the soliton equations along the principal orbits.

In \S 3 we focus, as mentioned above,
 on the case when the hypersurface is a circle bundle
over a product of (possibly inhomogeneous) K\"ahler-Einstein manifolds.
We find that the Ricci soliton equations, like the Einstein equations,
admit a class of explicit solutions representing solitons with K\"ahler metrics.
By judicious choice of parameters we can arrange the boundary conditions
so as to obtain steady and expanding solitons on vector bundles over products
of Fano K\"ahler-Einstein spaces (cf Thm \ref{completeGKRS}). These vector
bundles may be of rank one or of higher rank. The latter examples thus
generalise the $U(n)$-invariant K\"ahler-Ricci solitons on ${\mathbb C}^n$ due
to H. D. Cao \cite{Ca1}. However, the sectional curvatures in the new
examples are no longer positive.

We also find compact shrinking solitons, where the manifold is
a $\C \PP^1$-bundle over a product of Fano K\"ahler-Einstein manifolds, or
is obtained from such a bundle by blow-downs (cf Thm \ref{koiso}).
Furthermore, we produce examples of complete {\em noncompact} shrinking
solitons generalising those of \cite{FIK} (cf Thm \ref{FIK}).

In \S 4 we return to the strict cohomogeneity one setting and
consider principal orbits which are circle bundles over a generalised flag
variety. Here, as in \cite{DW}, we make the assumption that the isotropy
representation of the principal orbit is multiplicity free. We note that
the hypersurface is now a Riemannian submersion over a K\"ahler (though not
necessarily Einstein) metric.  We shall show that the soliton equations are
identical to those in \S 3, and hence we obtain new examples of steady,
expanding, and shrinking solitons in this setting also.

\section{\bf The Cohomogeneity One Ricci Soliton Equations}

In this section we will adapt the Ricci soliton equation (\ref{formRS})
to the cohomogeneity one setting, following basically the approach and
notation of \cite{EW}.

Accordingly, let $G$ be a compact Lie group acting via isometries on an
$n+1$-dimensional connected Riemannian manifold $(\overline{M}, \bar{g})$ with
one-dimensional orbit space, which is further assumed to be not a circle.
We choose a unit speed geodesic $\gamma(t)$ which intersects all principal
orbits orthogonally. Let $K$ denote the principal isotropy group
along $\gamma(t)$. Then there is an equivariant diffeomorphism
$$ \Phi: I \times (G/K) \longrightarrow M_0 $$
given by $\Phi(t, gK) = g\cdot \gamma(t)$, where $M_0 \subset \overline{M}$
is the open and dense subset consisting of all points lying on principal orbits
and $I$ is an open interval. We denote by $P_t$ the principal orbit passing
through $\gamma(t)$ and by $P$ an abstract copy of the homogeneous space $G/K$,
which is assumed to be connected. It then follows that
\begin{equation} \label{metricform}
 \Phi^*(\bar{g}) = dt^2 + g_t
\end{equation}
where $g_t$ is a one-parameter family of $G$-invariant metrics on $P$. It will be
convenient to fix a background metric $b$ on $P$ which is induced by a bi-invariant
metric on $G$. We can then write
$$ g_t(X, Y) = b(q_t(X), Y), \ \  X, Y \in TP $$
where $q_t$ is a $b$-symmetric automorphism of $TP$.

We will adopt the convention that $\overline{\nabla}$, $\overline{\rm Ric},$
and $\overline{R}$ denote respectively the Levi-Civita connection, the
Ricci tensor, and the scalar curvature of the metric $\overline{g}$, while
$\nabla^t$,  ${\rm Ric}^t,$ and $R^t$ denote the corresponding objects for $g_t$.
Whenever the context is clear, we will drop the $t$-dependence to simplify
the notation. In this spirit, we let $L_t$ denote the shape operator of
the orbit $P_t$, i.e., for any $X \in TP_t$
$$ L_t(X) := \overline{\nabla}_X N$$
where $N= \Phi_*(\frac{\partial}{\partial t})$ is a unit $G$-equivariant
normal field along $P_t$ with $\overline{\nabla}_N N =0$. Via the
diffeomorphism $\Phi$, we can regard $L_t$ as a one-parameter family of
$G$-equivariant, $g_t$-symmetric endomorphisms of $TP$. In particular,
its trace ${\rm tr}(L)$ is constant along $P_t$. We also have, for $X, Y \in TP$,
$$ \dot{g}_t(X, Y) = 2g_t(L_t(X), Y) $$
where $\ \dot{} \ $ denotes $\frac{d}{dt}$, and
$$ \Phi_* \dot{L} = \overline{\nabla}_N L. $$

In \cite{EW}, viewing $M_0$ as an equidistant hypersurface family, and using the
Gauss and Codazzi equations together with the Riccati equation for $L$, one
obtains
\begin{eqnarray}
\overline{\rm Ric}(N, N) & = & -{\rm tr}(\dot{L}) - {\rm tr}(L_t^2) \label{ricNN} \\
\overline{\rm Ric}(X, N) & = & -g_t(\delta^{\nabla^t} L_t, X) - d({\rm tr}L_t)(X) \label{ricXN}  \\
\overline{\rm Ric}(X, Y) & = & {\rm Ric}_t (X, Y) - {\rm tr}(L_t) g_t(L_t(X), Y) - g_t(\dot{L}(X), Y) \label{ricXY}.
\end{eqnarray}
where $X, Y \in TP_t$, $L_t$ is viewed as a $TP$-valued $1$-form on $TP$, and
$\delta^{\nabla^t}: T^*(P) \otimes TP \rightarrow TP$ is the codifferential.

Let us now consider the Ricci soliton equation on $(\overline{M}, \overline{g})$.
This becomes
\begin{equation} \label{RS}
\overline{\rm Ric}(\overline{g}) + \overline{\delta}^*{\overline{\omega}}
          + \frac{\epsilon}{2} \ \overline{g} = 0
\end{equation}
where $\overline{\omega}$ is a $1$-form on $\overline{M}$ and
$$\overline{\delta}^*: \Omega^1(\overline{M}) \longrightarrow S^2(T^*(\overline{M}))$$
is the symmetrized covariant differential. Note that if we take the vector field
$\overline{g}$-dual to $\overline{\omega}$ and use Lemma 1.60 in \cite{Be}
we obtain Eq.(1.8) on p. 4 of \cite{Cetc}. Note also that if we add a $1$-form
corresponding to a Killing vector field to $\overline{\omega}$ we obtain another solution
of the soliton equation.

Suppose next that $G$ is any compact group of isometries of $\overline{g}$.
We fix a unit volume Haar measure $d\mu_G$ on $G$. For $\varphi \in G$,
since $\varphi^* \overline{\delta}^* \overline{\omega} = \overline{\delta}^* \varphi^* \overline{\omega},$
we obtain from Eq.(\ref{RS}) that
$$\overline{\delta}^* \varphi^* \overline{\omega} = \overline{\delta}^* \overline{\omega},$$
so that $(\overline{M}, \overline{g}, \varphi^*{\overline{\omega}})$
is also a Ricci soliton. Moreover, for $p \in \overline{M}, X \in T_p{\overline{M}}$, we
may define
$$\tilde{\omega}_p(X) := \int_G \overline{\omega}_p(\varphi_* X)\ d\mu_G(\varphi). $$
By differentiating under the integral sign and the above fact, we readily obtain
$\overline{\delta}^*{\tilde{\omega}} =\overline{\delta}^* \overline{\omega}$.

In other words, if $(\overline{M}, \overline{g}, \overline{\omega})$ is a Ricci
soliton and $\overline{g}$ has a compact subgroup $G$ of isometries, we may assume
that the $1$-form $\overline{\omega}$ is also $G$-invariant. In the special case
of a gradient Ricci soliton with $\overline{\omega} = du$ for some smooth function
$u$ on $\overline{M}$, the above argument shows that we may assume that $u$ is
$G$-invariant. ($u$ is often called a {\em potential} for $\overline{\omega}$.)

Returning to the cohomogeneity one situation, we will assume from now on that
$\overline{\omega}$ is $G$-invariant. Then it follows that
\begin{equation} \label{omegaform}
\Phi^* \overline{\omega} = \xi(t)\ dt + \omega_t
\end{equation}
where $\xi$ is a function of $t$ only and $\omega_t$ is a $1$-parameter family of
$G$-invariant $1$-forms on $P$.

\begin{lemma} \label{diffomega} For $X, Y \in TP_t$, we have
\begin{enumerate}
  \item[a.] $(\overline{\delta}^* \overline{\omega})(N, N) = \dot{\xi},$
  \item[b.] $(\overline{\delta}^* \overline{\omega})(N, X) = \frac{1}{2} \ \dot{\omega}_t(X)
           - \omega_t(L_t(X)),$
  \item[c.] $(\overline{\delta}^* \overline{\omega})(X, Y) = \xi \ g_t(L_t(X), Y) +
                 (\delta^* \omega_t)(X, Y),$
\end{enumerate}
\end{lemma}

\noindent{\bf Proof.} Since $\overline{\nabla}_N N = 0$, it follows that
$$ (\overline{\delta}^* \overline{\omega})(N, N)  =  (\overline{\nabla}_N \overline{\omega})(N)
  = N(\overline{\omega}(N)) = \dot{\xi},$$
which gives the first assertion.

Next we consider $X \in TP$ and extend it first to a local vector field in $P$
and then via $\Phi$ to a local vector field in $M_0$. It follows that $[N, X] = 0$. Then we have
\begin{eqnarray*}
2(\overline{\delta}^* \overline{\omega})(N, X) & = & (\overline{\nabla}_N \overline{\omega})(X)
      + (\overline{\nabla}_X \overline{\omega})(N) \\
    & = & N({\omega}_t(X)) - \overline{\omega}(\overline{\nabla}_N X) + X(\overline{\omega}(N))
          -\overline{\omega}(\overline{\nabla}_X N)  \\
    & = & \dot{\omega}(X) -2 \overline{\omega}(\overline{\nabla}_X N)
\end{eqnarray*}
since $X(\xi) =0.$ This gives the second assertion.

For the third assertion, let us extend $X, Y$ to local vector fields as above.
We have
$$ \overline{\nabla}_X Y = \nabla_X Y + \ \overline{g}(\overline{\nabla}_X Y, N) N
      = \nabla_X Y - \ g(L(X), Y) N.$$
It follows that
\begin{eqnarray*}
2 (\overline{\delta}^* \overline{\omega})(X, Y)  & = & X(\overline{\omega}(Y)) + Y(\overline{\omega}(X))
       -\overline{\omega}(\overline{\nabla}_X Y + \overline{\nabla}_Y X) \\
   & = & X({\omega}(Y)) + Y({\omega}(X)) + 2\xi g(L(X), Y)
        - \omega(\nabla_X Y + \nabla_Y X) \\
    & = & 2 (\delta^* \omega)(X, Y) + 2\xi \ g(L(X), Y). \ \ \ \ \ \ \qed
\end{eqnarray*}

Combining the above Lemma with (\ref{ricNN})-(\ref{ricXY}), we obtain

\begin{prop} \label{RSeqn}
Let $(\overline{M}^{n+1}, \overline{g})$ be a Riemannian manifold
where $\overline{M}$ admits a cohomogeneity one action with respect to
some compact group $G$ of isometries of $\overline{g}$. Assume that
$\overline{\omega}$ is a $G$-invariant $1$-form on $\overline{M}$. Under the
diffeomorphism $\Phi : I \times P \rightarrow M_0$ induced by a unit-speed
geodesic $\gamma$ which intersects all principal orbits orthogonally, the
Ricci soliton equation for $\overline{g}$ and the vector field dual to
$\overline{\omega}$ is equivalent to
\begin{eqnarray}
-({\delta}^{{\nabla}^t} L_t)^{\flat}- d({\rm tr} (L_t))
      + \frac{1}{2} \ \dot{\omega_t} - \omega_t \circ L_t & = & 0, \label{RS-XN} \\
-{\rm tr}(\dot{L}) - {\rm tr}(L_t^2) + \dot{\xi} + \frac{\epsilon}{2} & = & 0, \label{RS-NN} \\
{\rm Ric}_t(X, Y) -\ {\rm tr}(L_t)\ g_t(L_t(X), Y) -g_t(\dot{L}(X), Y) \ + &  &  \label{RS-XY}  \\
     \xi \ g_t(L_t(X),Y) + (\delta^* \omega_t)(X, Y) + \frac{\epsilon}{2}\ g_t(X, Y) & = & 0, \nonumber
\end{eqnarray}
for all $X, Y \in TP_t, t \in I$.

Conversely, if $g_t$ $($resp. $\omega_t$$)$ is a smooth $1$-parameter family of
$G$-invariant metrics $($resp. $1$-forms$)$ on $P=G/K$ and $\xi$ is a smooth
function of $t\in I$ such that, with $L_t$ defined by $\dot{g}(X, Y) = 2 g_t(L_t(X), Y)$
for $X, Y \in TP$, the above system holds, then $\bar{g} = dt^2 + g_t$ and
$\bar{\omega} = \xi dt + \omega_t$ give a local Ricci soliton on $I \times P$.
\ \ \ \ \qed
\end{prop}

In the above we have used the notation $Z^{\flat}$ to denote the $1$-form
dual to the given vector field $Z$.

In the situation of a gradient Ricci soliton, $\overline{\omega} = du,$
where we may assume that $u$ is a $G$-invariant function,
Lemma \ref{diffomega} then implies that the only non-trivial components of
$\overline{\delta}^* \overline{\omega} = {\rm Hess}(u)$
are
\begin{eqnarray}
(\overline{\delta}^* \overline{\omega})(N, N) & = & \ddot{u} \\
(\overline{\delta}^* \overline{\omega})(X, Y) & = & \dot{u}\ g_t(L_t(X), Y)
\end{eqnarray}
where $X, Y \in TP_t$. We then obtain the system
\begin{eqnarray}
-({\delta}^{{\nabla}^t} L_t)^{\flat}- d({\rm tr} L_t) & = & 0, \label{GRS-XN} \\
-{\rm tr}(\dot{L}_t) - {\rm tr}(L_t^2) + \ddot{u} + \frac{\epsilon}{2} & = & 0, \label{GRS-NN} \\
{\rm Ric}_t(X, Y) -\ {\rm tr}(L_t)\ g_t(L_t(X), Y) -g_t(\dot{L}(X), Y) \ + &  &  \label{GRS-XY}  \\
     \dot{u} \ g_t(L_t(X),Y) +  \frac{\epsilon}{2}\ g_t(X, Y) & = & 0. \nonumber
\end{eqnarray}

\begin{rmk} \label{gen}
The above formulas are valid not just in the cohomogeneity one setting but
also when we have a manifold which, after removing some higher codimension
submanifolds, is an equidistant hypersurface family. More precisely, suppose
$\overline{M}^{n+1}$ has a smooth real-valued function $t$ with range an interval
$I$ such that $M_0 := t^{-1}({\rm int}(I))$ is diffeomorphic to a product
${\rm int}(I) \times P$ for a fixed $n$-dimensional manifold $P$. Furthermore,
assume that on $\overline{M}$ there is a Riemannian metric $\bar{g}$ and a $1$-form
$\overline{\omega}$ such that under pull-back via some such diffeomorphism they have
the respective forms (\ref{metricform}) and (\ref{omegaform}). Then
(\ref{ricNN})-(\ref{ricXY}) and Lemma \ref{diffomega} remain valid, and hence
so does Proposition \ref{RSeqn}.

Unlike the cohomogeneity one situation, however, $\tr(L_t)$ does not in
general depend only on $t$, so the $d(\tr(L_t))$ term in Eq.(\ref{RS-XN})
could be nonzero. As well, the scalar curvature of each level set $(P_t, g_t)$
is not necessarily constant. In the case of a gradient Ricci soliton,
(\ref{GRS-XN})-(\ref{GRS-XY}) require the assumption that the potential
$u$ is constant on each $P_t$.
\end{rmk}

\medskip

We now give a more precise description of the $G$-invariant $1$-forms
$\omega_t$, and hence of Eq.(\ref{RS-XY}). To this end, let
$$ \g = \kf \oplus \p $$
be an ${\rm Ad}_K$-invariant decomposition of $\g$ with respect to
the background bi-invariant metric on $G$, so that $\p \approx T_{[K]}(P)$.
Let $\p_0$ denote the subspace of $\p$ on which ${\rm Ad}_K$ acts as the identity.
Then $\omega_t$ can be regarded as a path in $\p_0^*$.

Let $\{Z_1, \cdots, Z_{\ell} \}$ be a $b$-orthonormal basis of $\p_0$ and
$\{\theta_1, \cdots, \theta_{\ell} \}$ be the dual basis. Then
$ \omega_t = \sum_{i=1}^{\ell} w_i(t) \theta_i$ for certain smooth
functions $w_i(t)$.
We will write down Eq.(\ref{RS-XY}) along the geodesic $\gamma(t)$.
We begin with the proof of part c of Lemma \ref{diffomega} where we
now assume that $X, Y$ are the Killing vector fields on $M_0$ generated by
two vectors in $\p \approx T_{[K]}(G/K)$. Using the Lie derivative,
we have
$$ 2 \overline{\delta}^* \overline{\omega}(X, Y) = 2 \xi \ g_t(L_t(X), Y) +
      ({\mathscr L}_X \omega_t)(Y) + ({\mathscr L}_Y \omega_t)(X) -
       \omega_t(\nabla^t_Y X + \nabla^t_X Y). $$
By $G$-equivariance, the second and third terms on the right vanish,
while the last term can be written as
$$-\sum_i \ b(\nabla^t_X Y + \nabla^t_Y X, Z_i)\ w_i(t)
    = - g_t(\nabla^t_X Y + \nabla^t_Y X, W_t) $$
where $W_t:= q_t^{-1}(\sum_i w_i(t)Z_i)$ is the $g_t$-dual to $\omega_t$.
By Proposition 7.28 in \cite{Be}, the last expression becomes
$$ - g_t([W_t, X]_{\p}, \ Y) - g_t(X, \ [W_t, Y]_{\p}), $$
where we have switched to using brackets in the Lie algebra $\g$.
Hence Eq.(\ref{RS-XY}) becomes
\begin{eqnarray}
{\rm Ric}_t(X, Y) + (\xi -{\rm tr}(L))\ g_t(L(X), Y) - g_t(\dot{L}(X), Y)
          + \frac{\epsilon}{2} \ g_t(X, Y)  &  & \label{RS-XY2} \\
    - \frac{1}{2}\left(g_t([W_t, \ X]_{\p}, Y) + g_t(X,\ [W_t, Y]_{\p})\right)
      & =& 0 \nonumber.
\end{eqnarray}

\begin{rmk}
Note that in the last two terms of (\ref{RS-XY2}) we have the
operator $\frac{1}{2}({\rm ad}_{W_t} + {\rm ad}^*_{W_t})$
where ${}^*$ denotes the adjoint with respect to $g_t$. To deal with it,
the following observations are useful.
\begin{enumerate}
\item[(i)] Using the background biinvariant metric $b$ we can decompose
$\p$ as an ${\rm Ad}_K$-invariant orthogonal sum
$$ \p = \p_0 \oplus \p_1 \oplus \cdots \oplus \p_r $$
where $\p_i$ and $\p_j$ are inequivalent  orthogonal real representations
if $i\neq j$ and each $\p_i$ is a sum of equivalent irreducible orthogonal summands.
Then $\p_0$ is an ${\rm Ad}_K$-invariant subalgebra of $\g$. In fact it is
the Lie algebra of $C^0(K)/Z^0(K)$, where $C$ and $Z$ denote respectively
the centralizer and centre, so $\p_0$ is of compact type.
\item[(ii)] We have $[\p_0, \p_i] \subset \p_i$ for all $i$ because
for any irreducible ${\rm Ad}_K$-submodule $\m$ of $\p$ and any $h \in C(K)$,
${\rm Ad}_h$ induces an ${\rm Ad}_K$-equivariant isomorphism of $\m$
with ${\rm Ad}_h(\m)$, so the latter lies in the same $\p_j$ as $\m$ does.
Note that $q_t(\p_i) \subset \p_i$ for all $i$ as well. Since $W_t \in \p_0$,
it follows that the operator $\frac{1}{2}({\rm ad}_{W_t} + {\rm ad}^*_{W_t})$
maps $\p_i$ to itself.
\end{enumerate}
\end{rmk}

\begin{rmk} \label{commondecomp}
In studying cohomogeneity one metrics one sometimes makes the assumption
that there is a $b$-orthogonal decomposition of $\p$ into ${\rm Ad}_K$-invariant
summands, say $\p = \m_1 \oplus \cdots \oplus \m_{\ell}$, such that
$g_t |_{\m_i} = f_i(t)^2 \, b|_{\m_i}$. (Here one does not assume that $\m_i$ is
irreducible.) This is the situation in a multiple warped product or
when $\p$ is a sum of pairwise inequivalent irreducible orthogonal
representations. Suppose further that ${\rm ad}_{W_t}$ preserves the above
decomposition.  Then   ${\rm Ad}_{\exp(sW_t)}$ leaves all $\m_i$ invariant,
and so it is clear from the form of the metric that $\exp(sW_t)$ acts on the right
of $G/K$ via isometries of $g_t$. Hence
$\delta^* \omega_t = \frac{1}{2}{\mathscr L}_{W_t} g_t = 0$.
Note that by Prop. 3.18 in \cite{BB}, we also have $\overline{\ric}(X, N) = 0.$
So by (\ref{RS-XN}), (\ref{ricXN}), and (\ref{diffomega}) we are reduced to
the case of a gradient Ricci soliton.
\end{rmk}

\medskip Next we consider briefly the smoothness criterion for
$\overline{\omega}$ when there is a special orbit $G/H$, where $K
\subset H$ and $H/K = S^k$.  This is obtained in essentially the same
way as for the metric $\overline{g}$ (cf. the first remark on p. 114
and section 1 in \cite{EW}).  Let $V \approx \R^{k+1}$ denote a normal
slice at $[H] \in G/H$ and $\langle \ , \ \rangle_0$ denote the
Euclidean metric on it (induced by $\overline{g}$ if it is already
given). The subgroup $H$ acts orthogonally irreducibly on $V$ and,
because of the cohomogeneity one assumption, transitively on the unit
sphere in it. But in general $H$ need not act effectively on $V$. We
denote by $\p_{+}$ and $\p_{-}$ respectively the subspaces of the
tangent space $T_{eK}(G/K)$ corresponding to $H/K$ and $G/H$.

On the tube around $G/H$, the $1$-form $\overline{\omega}$ is determined
by an $H$-equivariant map from $V  \rightarrow V^* \oplus \ \p_{-}^*$ and vice
versa. The Taylor series for this map gives rise to $H$-equivariant homogeneous
polynomials on $V$ with values in $V^* \oplus \ \p_{-}^*$, i.e., elements
of $\mathscr{W}_m:={\rm Hom}_H(S^m(V), V^* \oplus \p_{-}^*)$. For $\overline{\omega}$ of the
form $\xi(t) dt + \omega_t$ to be smooth, it is necessary and sufficient
(cf Lemma 1.1 in \cite{EW}) that, for all $m$,  its $m$th Taylor coefficient
(as a function of $t$) be the restriction to the unit sphere in $V$ of elements
of $\mathscr{W}_m$. Note that at $[K] \in H/K =S^{k}$ these elements
must have a value which is fixed by $K$.

Besides determining the spaces $\mathscr{W}_m$, one of the technical points
of applying the above smoothness criterion is to reconcile the relationship
between the biinvariant metric $b$, the metric $dt^2 + g_t$, and the
Euclidean metric $\langle \ , \  \rangle_0$ on $V$. It is now convenient to
write
$$ \g = \kf \oplus \p_{+} \oplus \p_{-} $$
where the above decomposition is $b$-orthogonal and $\h=\kf \oplus \p_{+}$.
We have a corresponding orthogonal decomposition $\p_0 = \af_{+} \oplus \af_{-}$,
and we will take the basis $\{Z_1, \cdots, Z_{\ell}\}$ of $\p_0$ chosen before
to be adapted to this decomposition.

Note that if $\p_0 = 0$, then on $I \times P$ we have $\overline{\omega} = du$,
where $u(t)$ is an anti-derivative of $\xi(t)$.
Furthermore, if there is a special orbit $G/H$ (with $K \subset H$) then
smoothness further implies that  $u$ is even in $t$.
We therefore obtain

\begin{lemma} \label{gradRS}
Suppose that $(\overline{M}, \overline{g}, \overline{\omega})$
is a cohomogeneity one Ricci soliton with $\overline{\omega}$ chosen to be
$G$-invariant. If the isotropy representation of the principal orbit has no
trivial irreducible summands, then it must be a gradient Ricci soliton.

For a gradient Ricci soliton of cohomogeneity one with $G$-invariant potential,
we have $\overline{\omega} = \dot{u} \ dt$ for some smooth function $u(t)$.
If there is further a special orbit at $t=0$, then
$u(t)$ must be even in $t$.   \ \ \ \qed
\end{lemma}

\section{\bf Some Consequences of the Contracted Second Bianchi Identity}

We begin with a general remark about the Ricci soliton equation (\ref{formRS}).
It is well-known (see, e.g., Lemma 1.10 in \cite{Cetc}) that the contracted second
Bianchi identity yields the consequence
\begin{equation} \label{elliptic}
\Delta \omega = 2 \omega \circ r
\end{equation}
where $\Delta$ is the Laplace-Beltrami operator and $r$ is the Ricci endomorphism.
The following result is a simple extension of Theorem 5.1 in \cite{DTK} and
a special case of Bando's result \cite{Ban} for the Ricci flow.

\begin{lemma} \label{DTKz}
Assume that ${\rm dim}\, M \geq 3$. Then the system consisting of $($\ref{formRS}$)$
and $($\ref{elliptic}$)$ is a quasi-linear elliptic system in harmonic coordinates.
Hence a $C^{2}$ solution $(g, \omega)$ of the Ricci soliton equation is
automatically real analytic.
\end{lemma}

\noindent{\bf Proof.}  Note that in harmonic coordinates, the principal symbol
of the linearization of the system at a solution  is given by
$$ (\rho, \eta) \mapsto (\frac{1}{2}|\zeta|_g^2 \, \rho, \, \frac{1}{2}|\zeta|_g^2
\, \eta + \cdots )$$
where $(\rho, \eta) \in S^2T^*M \oplus T^*M,$  $\cdots$ denotes an expression linear
in $\rho$,  and $\zeta$ is a nonzero cotangent vector. Hence the system is
quasi-linear elliptic. If $(g, \omega)$ are $C^2$ in harmonic coordinates,
we can apply Morrey's interior regularity theorem as in \cite{DTK}. However,
in transforming $(g, \omega)$ to harmonic coordinates, the resulting tensors
are only in $C^{1, \alpha}$. We then need to apply Theorems 8.8 and 9.19 in
\cite{GT} to the components of the system to see that in fact $(g, \omega)$
are $C^{2, \alpha}$.  \qed

\medskip
Returning to the cohomogeneity one situation, if we use the connection Laplacian
instead of the Laplace-Beltrami operator, (\ref{elliptic}) becomes
\begin{equation} \label{oneform}
 \overline{\nabla}^* \overline{\nabla} \overline{\omega} = \overline{\omega}\circ \bar{r}
\end{equation}
where $\bar{r}$ is now the Ricci endomorphism corresponding to the Ricci
tensor $\overline{\rm Ric}$ via $\bar{g}$.  It is classical \cite{Ya} that this
equation is also satisfied for a $1$-form dual to a Killing field.  So in considering
symmetric Ricci solitons, in addition to its role in ellipticity, it is natural
to examine the relation of (\ref{oneform}) to the Ricci soliton equations.

In this section we show that Eqs (\ref{RS-XN}) and (\ref{RS-NN}) can be
replaced by Eq.(\ref{oneform}), provided there is a special orbit of dimension
strictly smaller than that of a principal orbit and provided $C^{2}$
regularity for $\bar{g}$ and $\overline{\omega}$ has been established.
This generalizes an observation of A. Back for the cohomogeneity
one Einstein equations (cf \cite{Ba} and \cite{EW}, pp. 118-120).

\begin{lemma} \label{oneformparts}
For a $G$-equivariant $1$-form $\overline{\omega}$ of the form $($\ref{omegaform}$)$,
Eq.$($\ref{oneform}$)$ is equivalent to the equations
\begin{eqnarray}
2 \sum_i (\delta^* \omega_t)(L_t(e_i), e_i)  &=& \ddot{\xi} + {\rm tr}(L_t) \dot{\xi}
     -\left({\rm tr}(\dot{L}) + 2\ {\rm tr}(L_t^2)\right) \xi, \label{oneform-N} \\
\nabla^* \nabla \omega_t - \omega_t \circ r_t & = & \ddot{\omega}-2 (\omega_t \circ L_t)^{\cdot}
    + {\rm tr}(L_t) \dot{\omega}  -2{\rm tr}(L_t)(\omega_t \circ L_t)
    -2 \xi (\delta^{\nabla^t}L_t)^{\flat}  \label{oneform-X}
\end{eqnarray}
where $\{e_1, \cdots, e_n\}$ is an orthonormal basis of $TP_t$.
\end{lemma}

\noindent{\bf Proof.} This proceeds essentially by straight-forward computation
using a local orthonormal frame $\{N, e_1, \cdots, e_n\}.$ First, we consider
\begin{eqnarray*}
\overline{\nabla}^* \overline{\nabla} \overline{\omega}(N) &=& - (\overline{\nabla}_N
     \overline{\nabla}_N \overline{\omega})(N)
     - \sum_i \left((\overline{\nabla}_{e_i} \overline{\nabla}_{e_i} \overline{\omega})(N)
       -(\overline{\nabla}_{\overline{\nabla}_{e_i} e_i} \overline{\omega})(N) \right) \\
  &=& -N(N\xi)
   +\sum_i \left( e_i ({\omega}_t(L_t(e_i)))+ ({\overline{\nabla}}_{e_i} \overline{\omega})(L_t(e_i))
     -\omega_t(L_t(\nabla^t_{e_i} e_i)) + \bar{g}({\overline{\nabla}}_{e_i} e_i, N) \dot{\xi} \right) \\
  &=& -\ddot{\xi} + \sum_i \left(2({\nabla}_{e_i} \omega_t)(L_t(e_i))+ \omega_t((\nabla_{e_i} L)(e_i))
        - \bar{g}(\overline{\nabla}_{e_i}(L_t(e_i)), N) \xi \right) -{\rm tr}(L)\dot{\xi} \\
  &=& -\ddot{\xi} +2 \sum_{i,j} g_t(L_t(e_i), e_j)(\nabla_{e_i}^t \omega_t)(e_j)
        -\omega_t(\delta^{{\nabla}^t} L_t) +\sum_i g_t(L_t(e_i), L_t(e_i)) \xi -{\rm tr}(L) \dot{\xi} \\
  &=& -\ddot{\xi} + 2 \sum_i (\delta^* \omega_t)(L_t(e_i), e_i) - \omega_t(\delta^{{\nabla}^t} L_t)
        + {\rm tr}(L_t^2) \xi - {\rm tr}(L_t) \dot{\xi}. \\
\end{eqnarray*}
Also, using (\ref{ricNN}) and (\ref{ricXN}), we have
\begin{eqnarray*}
\overline{\omega}(\bar{r}(N)) &=&  \overline{\rm Ric}(N, N)\xi  +
        \sum_i \ \overline{\rm Ric}(N, e_i) \ \omega_t(e_i) \\
     &=& -({\rm tr}(\dot{L}) + {\rm tr}(L_t^2))\xi
          -\sum_i \omega_t(e_i)\left(g_t(\delta^{{\nabla}^t}L_t, e_i) + e_i({\rm tr}(L_t))  \right).
\end{eqnarray*}
Putting the two computations together and noting that ${\rm tr}(L_t)$ depends only on
$t$, we obtain (\ref{oneform-N}).

Similarly, by using (\ref{ricXN}) and (\ref{ricXY}) we have, for $X \in TP_t$ (extended
to a local vector field in the usual manner such that $[X, N] =0$),
\begin{eqnarray*}
\overline{\omega}(\bar{r}(X)) &=& \overline{\rm Ric}(X, N) \xi +
            \sum_i \ \overline{\ric}(X, e_i)\ \omega_t(e_i)        \\
   &=& -g_t(\delta^{{\nabla}^t} L_t, X) \xi + \sum_i \omega_t(e_i) \left(\ric(X, e_i)
            - {\rm tr}(L_t)g_t(L_t(X), e_i) - g_t(\dot{L}(X), e_i) \right) \\
   &=& -g_t(\delta^{{\nabla}^t} L_t, X) \xi + \omega_t (r_t(X))- \tr (L_t)\ \omega_t(L_t(X))
           -\omega_t(\dot{L}(X)).
\end{eqnarray*}
On the other hand, we have
\begin{eqnarray*}
(\overline{\nabla}^* \overline{\nabla} \overline{\omega})(X) &=&
     - (\overline{\nabla}_N \overline{\nabla}_N \overline{\omega})(X)
     - \sum_i \left((\overline{\nabla}_{e_i} \overline{\nabla}_{e_i} \overline{\omega})(X)
       -(\overline{\nabla}_{\overline{\nabla}_{e_i} e_i} \overline{\omega})(X) \right) \\
  & =& -N((\overline{\nabla}_N \overline{\omega})(X))
          +(\overline{\nabla}_N \overline{\omega})(L_t(X))  \\
  &  &   -\sum_i \left( e_i((\overline{\nabla}_{e_i} \overline{\omega})(X))
    -(\overline{\nabla}_{e_i} \overline{\omega})(\overline{\nabla}_{e_i} X)
    -(\overline{\nabla}_{e_i}e_i)(\omega_t(X))
    +\overline{\omega}(\overline{\nabla}_{\overline{\nabla}_{e_i} e_i} X) \right).
\end{eqnarray*}

After systematically splitting the covariant derivatives in the above expression
into their components along and orthogonal to the principal orbits, we obtain
$$ -\ddot{\omega}(X) + 2N(\omega_t (L_t(X)) -\omega_t(\dot{L}(X))
  -\sum_i (\nabla_{e_i}^t \nabla_{e_i}^t \omega_t)(X) - \xi \sum_i \left(e_i(g_t(X, L_t(e_i)))
     + \bar{g}(\overline{\nabla}_{e_i}(\nabla_{e_i}^t X), N)\right)  $$
$$   -\tr(L_t) \dot{\omega}(X)  +\tr(L_t)\omega_t(L_t(X)) + \sum_i (\nabla_{e_i}^t e_i)(\omega_t(X))
    + \sum_i \bar{g}(X, L_t(\nabla_{e_i}^t e_i)) \xi -\sum_i \omega_t(\nabla^t_{\nabla_{e_i}^t e_i} X). $$
Unravelling the covariant derivatives further in the above, we arrive at
$$ -\ddot{\omega}(X) + 2 (\omega_t \circ L_t)^{\cdot}(X)- \omega_t(\dot{L}(X))
    + (\nabla^* \nabla \omega_t)(X)  +\xi g_t(X, \delta^{{\nabla}^t}L_t)
  + \tr(L_t) \omega_t(L_t(X)) - \tr(L_t)\dot{\omega}(X). $$
Combining this with the computation for $\overline{\omega}(\bar{r}(X))$ above, we
finally obtain (\ref{oneform-X}). \ \ \ \ \qed

\begin{rmk} Notice that in the setting of Remark \ref{gen}, the above
Lemma remains valid provided that all the hypersurfaces $P_t$ have constant
mean curvature.
\end{rmk}

\begin{rmk} \label{cor-1formparts}
In the case of a gradient Ricci soliton with $\overline{\omega} = \dot{u} \ dt$,
Eq.$($\ref{oneform}$)$ becomes the system
\begin{eqnarray}
 \frac{d^3 u}{dt^3} + \tr(L_t) \ \ddot{u}
     -(\tr (\dot{L}) + 2 \ \tr(L_t^2))\ \dot{u} &=& 0  \label{1formGRS-NN} \\
  \dot{u} \ g_t(\delta^{{\nabla}^t} L_t, X) &=& 0 \label{1formGRS-XN}
\end{eqnarray}
for all $X \in TP_t $ and all $t$.
In particular, for a non-trivial gradient Ricci soliton, the second equation
coincides with Eq.$($\ref{GRS-XN}$)$ $($since $\tr(L_t)$ depends only on $t$$)$.
If we combine Eq.(\ref{1formGRS-NN}) with Eq.(\ref{GRS-NN}) then we obtain
\begin{equation} \label{altGRS-NN}
\frac{d^3 u}{dt^3} + \tr(L_t) \ \ddot{u}
     +\tr (\dot{L})\, \dot{u}  - 2 \ddot{u} \dot{u} -\epsilon \dot{u}  = 0.
\end{equation}
This last equation, which may be viewed as an equation for $\dot{u}$,
corresponds to the first integral observed in \cite{Iv} (p. 242) and more
generally in \cite{Ha} (pp. 84-85) and \cite{Ca3} (p. 123).
\end{rmk}

\begin{prop} \label{bianchi-RS}
Suppose that $\bar{g}$ and $\overline{\omega}$ are respectively a
$G$-invariant metric and  $1$-form on ${\rm int}(I) \times P$ of the form
$($\ref{metricform}$)$ and $($\ref{omegaform}$)$. Let $v$ denote the ratio of
the volume of $g_t$ to that of the background metric $b$.

\smallskip
{\noindent $($i$)$} If Eqs.$($\ref{RS-XY}$)$ and $($\ref{oneform-X}$)$ hold,
then
   $$ \frac{\partial}{\partial t}  \left( v \left(\overline{\ric}(N, X) +
        \frac{1}{2}\ \dot{\omega}(X) - \omega_t (L_t(X)) \right)\right) = 0 $$
along $\gamma(t)=(t, x_0)$, where $X \in T_{x_0}(P)$ and $x_0 \in P$ are arbitrary.

\smallskip
{\noindent $($ii$)$} If Eqs.$($\ref{RS-XN}$)$, $($\ref{RS-XY}$)$ and
      $($\ref{oneform-N}$)$ hold, then along each $\gamma(t)$ as above we have
$$ \frac{\partial}{\partial t} \left(v^2 \left(\overline{\ric}(N, N) + \dot{\xi}
      + \frac{\epsilon}{2} \right) \right)=0. $$
\end{prop}

\smallskip
{\noindent {\bf Proof of (i).}}
By $G$-equivariance, we may assume that $x_0$ is the base point $[K] \in P =G/K$.
We can extend $X$ to a vector field in $P$ near $x_0$ and then to a local
vector field in ${\rm int}(I) \times P$ in a neighbourhood of $\gamma$. When
we compute at $\gamma(t)$, for a fixed value of $t$,
we will choose a local $g_t$-orthonormal moving frame
$\{e_1, \cdots, e_n\}$ on $P$ such that at $x_0$ we have
$\nabla^t_{e_i} e_j = 0$ for all $1 \leq i, j \leq n$. This then induces a moving frame
on ${\rm int}(I) \times P$ which is not necessarily orthonormal off $P_t$ but which still
commutes with $N$.

Since the scalar curvature of $\bar{g}$ is constant along $P_t$, the contracted second
Bianchi identity yields at $\gamma(t)$
$$ 0 =(\overline{\nabla}_N \overline{\ric})(N, X) +
      \sum_i (\overline{\nabla}_{e_i} \overline{\ric})(e_i, X). $$
As $\overline{\nabla}_N N= 0$ and $[N, X]= 0$, the first term on the right becomes
$$ N(\overline{\ric}(N, X)) - \overline{\ric}(N, L_t(X)).$$
The second term equals
$$ \sum_i \left(e_i(\overline{\ric}(e_i, X)) - \overline{\ric}(e_i, \overline{\nabla}_{e_i} X)
          - \overline{\ric}(\overline{\nabla}_{e_i} e_i, X) \right). $$
We can now split the ambient covariant derivatives into their components along
and orthogonal to $P_t$, after which we may apply Eq.(\ref{RS-XY}) three times to get
\begin{eqnarray*}
 \overline{\ric}(L_t(X), N)) & + & \tr(L_t)\overline{\ric}(X, N)
     - \sum_i e_i \left(\xi g_t(L_t(X), e_i ) + (\delta^* \omega_t)(X, e_i)\right) \\
     & + &  \sum_i \left( \xi g_t(L_t(e_i), \nabla^t_{e_i} X) +
           (\delta^* \omega_t)(e_i, \nabla^t_{e_i} X) \right)
 \end{eqnarray*}

Combining the above computations and unravelling the differentiation with respect to $e_i$ yield
\begin{eqnarray*}
0 &= & N (\overline{\ric}(N, X)) + \tr(L_t)\ \overline{\ric}(X, N)
      - \sum_i \xi \ g_t((\nabla^t_{e_i}L)(X), e_i) \\
  &  & +\frac{1}{2} \sum_i \left( -(\nabla_{e_i}^t \nabla_X^t \omega_t)(e_i) -
         (\nabla^t_{e_i} \nabla^t_{e_i} \omega_t)(X) +
         (\nabla^t_{\nabla_{e_i}^t X} \omega_t) (e_i) \right)  \\
  & = & N( \overline{\ric}(N, X)) + \tr(L_t)\ \overline{\ric}(X, N)
         + \xi \ g_t(X, \delta^{{\nabla}^t} L_t) + \frac{1}{2} (\nabla^* \nabla \omega_t)(X) \\
   &  & + \frac{1}{2} \sum_i \left(-(\nabla^t_X \nabla^t_{e_i} \omega_t)(e_i)
        - (\nabla^t_{[e_i, X]} \omega_t)(e_i)  - \omega_t( R^t_{X, e_i} (e_i)))
           + (\nabla^t_{\nabla_{e_i}^t X} \omega_t) (e_i)\right)  \\
& = & N(\overline{\ric}(N, X)) + \tr(L_t)\ \overline{\ric}(X, N)
         + \xi \ g_t(X, \delta^{{\nabla}^t} L_t)  + \frac{1}{2} X(\delta^{{\nabla}^t} \omega_t) \\
  &  &  + \frac{1}{2} \left((\nabla^* \nabla \omega_t)(X)- \omega_t(r_t(X)) \right)
\end{eqnarray*}

Note that $\delta^{{\nabla}^t}\omega_t$ is a $G$-invariant function, so
$X(\delta^{{\nabla}^t} \omega_t)=0$. We now multiply the last equation by
$v$ and use the fact that $\dot{v} = \tr(L_t) v$ to get
\begin{equation} \label{bianchi-1}
0= \frac{\partial}{\partial t} \left(v( \overline{\ric}(X, N)) \right)
     + v \xi \ g_t(X, \delta^{{\nabla}^t} L_t)
     +   \frac{v}{2} \left((\nabla^* \nabla \omega_t)(X)- \omega_t(r_t(X)) \right).
\end{equation}
If (\ref{oneform-X}) holds, then we obtain
$$0= \frac{\partial}{\partial t} \left(v( \overline{\ric}(X, N)) \right)
    + \frac{\partial}{\partial t} \left( v \left( \frac{1}{2}\, \dot{\omega}(X)
       -  (\omega_t \circ L_t)(X) \right)\right), $$
as required. \ \ \ \ \qed

\medskip
{\noindent{\bf Proof of (ii).}} Recall that we assume (\ref{RS-XN}) and (\ref{RS-XY})
in the following. This time let $\{N, e_1, \cdots, e_n\}$ be an orthonormal frame in a
neighbourhood of $\gamma(t)$ adapted to ${\rm int}(I) \times P$. By the contracted
second Bianchi identity, we have
\begin{eqnarray*}
\frac{1}{2} \, d\bar{R}(N) & = & N(\overline{\ric}(N, N)) + \sum_i \left(e_i(\overline{\ric}(e_i, N))
      - \overline{\ric}(\overline{\nabla}_{e_i} e_i, N) - \overline{\ric}(e_i, L_t(e_i)) \right) \\
  & = & N(\overline{\ric}(N, N)) + \tr(L_t)\; \overline{\ric}(N, N) + \xi \, \tr(L_t^2)
           + \frac{\epsilon}{2} \tr(L_t) + \sum_i (\delta^* \omega_t)(L_t(e_i), e_i) \\
   &  &  + \sum_i \left(e_i(-\frac{1}{2}\dot{\omega}(e_i)) + \omega_t(L_t(e_i))) +
            \frac{1}{2}\dot{\omega}(\nabla^t_{e_i} e_i) - (\omega_t \circ L_t)(\nabla^t_{e_i} e_i) \right)
\end{eqnarray*}
where we have used (\ref{RS-XN}) twice and (\ref{RS-XY}) once. On the other hand,
\begin{eqnarray*}
 d\bar{R}(N) & = & N\left(\overline{\ric}(N, N) + \sum_i \overline{\ric}(e_i, e_i) \right)  \\
       & = & N(\overline{\ric}(N, N)) - \dot{\xi}\, \tr(L_t) - \xi\,  \tr(\dot{L}) -N(\tr(\delta^* \omega_t)),
\end{eqnarray*}
where we have used (\ref{RS-XY}) again.

The above expressions for $d\bar{R}(N)$ combine to give an equation, which,
after multiplying by $v^2$ can be written as follows:
\begin{eqnarray}
0 &=& N\left(v^2\left(\frac{\epsilon}{2} + \overline{\ric}(N, N)\right)\right)
           + v^2\left(2 \xi \, \tr(L_t^2) +(\xi \tr(L_t))^{\cdot} \right) \label{bianchi-2}\\
  &  & + v^2 \left( 2\sum_i (\delta^* \omega_t)(L_t(e_i), e_i)
         + N(\tr(\delta^* \omega_t)) + \delta^{{\nabla}^t} \dot{\omega}
         -2 \delta^{{\nabla}^t}(\omega_t \circ L_t)  \right). \nonumber
\end{eqnarray}

One readily checks that
\begin{eqnarray*}
     \delta^{{\nabla}^t}(\omega_t \circ L_t) & = & \omega_t(\delta^{{\nabla}^t} L_t)
       -\sum_i (\nabla^t_{e_i} \omega_t)(L_t(e_i))  \\
      & = & \omega_t(\delta^{{\nabla}^t} L_t) - \sum_i (\delta^* \omega_t)(L_t(e_i), e_i).
\end{eqnarray*}

In order to compute $N(\tr(\delta^* \omega_t)) + \delta^{{\nabla}^t} \dot{\omega}$,
we switch to a frame $\{e_1, \cdots, e_n \}$, orthonormal at $\gamma(t)$, of the type
used in the proof of (i) above. Then
$$ N(\tr(\delta^* \omega_t))= \frac{1}{2} N \left(\sum_{i, j} g_t^{ij}\left((\nabla^t_{e_i} \omega_t)(e_j)
       + (\nabla^t_{e_j} \omega_t)(e_i)\right)\right). $$
Note that at $\gamma(t)$, we have
$$\dot{g}^{ij} = - \dot{g}_{ij} = -2 g_t(L_t(e_i), e_j)$$
and
$$ \delta^{{\nabla}^t} \dot{\omega}= -\sum_i (\nabla^t_{e_i} \dot{\omega})(e_i)
       = -\sum_i (\delta^* \dot{\omega})(e_i, e_i) = -\tr(\delta^*{\dot{\omega}}).$$

Using $[e_i, N]= 0$ we now have
\begin{eqnarray*}
N((\nabla^t_{e_i} \omega_t)(e_j)) &=& N\left(e_i(\omega_t(e_j)) - \omega_t(\nabla^t_{e_i} e_j)\right) \\
    & = & e_i N(\omega_t(e_j)) -\dot{\omega}(\nabla^t_{e_i} e_j)
              -\omega_t(\overline{\nabla}_N(\nabla^t_{e_i} e_j)) \\
    & = & (\nabla^t_{e_i} \dot{\omega})(e_j)
           -\omega_t(\overline{\nabla}_N \overline{\nabla}_{e_i} e_j) \\
    & = & (\nabla^t_{e_i} \dot{\omega})(e_j) -\omega_t((\nabla^t_{e_i} L_t)(e_j)) -
              \omega_t(\overline{R}_{N, e_i} (e_j)).
\end{eqnarray*}
It follows that
\begin{eqnarray*}
 N(\tr(\delta^* \omega_t)) + \delta^{{\nabla}^t} \dot{\omega}
      & = & -\tr(\delta^* \dot{\omega}) -2 \sum_{i, j} g_t(L_t(e_i), e_j)\, (\delta^* \omega_t)(e_i, e_j) \\
    &  &  +\frac{1}{2} \sum_{i} \left( 2 (\delta^* \dot{\omega})(e_i, e_i)
           -2 \omega_t((\nabla^t_{e_i} L_t)(e_i)) - 2 \omega_t(\overline{R}_{N, e_i} (e_i))  \right) \\
    &=&  -2 \sum_i (\delta^* \omega_t)(L_t(e_i), e_i)  + \omega_t(\delta^{{\nabla}^t} L_t)
            - \sum_i \omega_t(\overline{\ric}(N, e_i) e_i) \\
    & = & -2 \sum_i (\delta^* \omega_t)(L_t(e_i), e_i) + 2 \omega_t(\delta^{{\nabla}^t} L_t)
\end{eqnarray*}
where we have used (\ref{ricXN}) in the last step.

Eq. (\ref{bianchi-2}) now becomes
$$ 0 = N\left(v^2\left(\frac{\epsilon}{2} + \overline{\ric}(N, N)\right)\right)
           + v^2\left(2 \xi \, \tr(L_t^2)
           +(\xi \tr(L_t))^{\cdot} + 2 \sum_i (\delta^* \omega_t)(L_t(e_i), e_i)\right).  $$
Applying (\ref{oneform-N}) we get
\begin{eqnarray*}
 0 &=& N\left(v^2\left(\frac{\epsilon}{2} + \overline{\ric}(N, N)\right)\right)
           + v^2\left(2 \xi \, \tr(L_t^2) +(\xi \tr(L_t))^{\cdot}
           + \ddot{\xi} + {\rm tr}(L_t) \dot{\xi}
     -\left({\rm tr}(\dot{L}) + 2\ {\rm tr}(L_t^2)\right) \xi \right) \\
   & = & N\left(v^2\left(\frac{\epsilon}{2} + \overline{\ric}(N, N)\right)\right)
          + v^2 (\ddot{\xi} + 2 {\rm tr}(L_t) \dot{\xi}) \\
    & = & N\left(v^2\left(\frac{\epsilon}{2} + \overline{\ric}(N, N)
                + \dot{\xi} \right)\right),
\end{eqnarray*}
as asserted. \ \ \ \ \ \qed

In either case of Proposition \ref{bianchi-RS}, if there is a special orbit
whose dimension is strictly smaller than that of a principal orbit, then
the volume distortion $v$ becomes zero at the special orbit. It follows that
if both the metric $\bar{g}$ and the $1$-form $\overline{\omega}$ is $C^2$,
then along $\gamma$, we have
$\overline{\ric}(N, X) + \frac{1}{2}\ \dot{\omega}(X) - \omega_t (L_t(X)) = 0$
(resp. $ \overline{\ric}(N, N) + \dot{\xi}+ \frac{\epsilon}{2} = 0 $).
By $G$-equivariance, these expressions are zero everywhere on $\overline{M}$.
So we have a $C^2$ solution of the elliptic system in Lemma \ref{DTKz}.
We have therefore deduced

\begin{cor} \label{cor-bianchi}
Let $G$ be a compact Lie group acting isometrically with cohomogeneity $1$
on a connected Riemannian manifold $(\overline{M}, \bar{g})$ where $\bar{g}$ is of
class $C^2$. Let $\overline{\omega}$ be a $G$-equivariant $1$-form of class
$C^2$. Suppose that on $M_0$ the equations $($\ref{RS-XY}$)$ and $($\ref{oneform-X}$)$
are satisfied. If $\overline{M}$ has a special orbit with dimension
strictly smaller than that of a principal orbit, then $($\ref{RS-XN}$)$
holds everywhere. If in addition $($\ref{oneform-N}$)$ holds on $M_0$,
then $($\ref{RS-NN}$)$ also holds everywhere.

In particular, if $($\ref{RS-XY}$)$, $($\ref{oneform-N}$)$, and $($\ref{oneform-X}$)$
hold, $(\overline{g}, \overline{\omega})$ is a real analytic Ricci soliton
on $\overline{M}$. \ \ \ \ \qed
\end{cor}

\begin{rmk} \label{gen2}
In the setting of Remark \ref{gen}, in order for Proposition \ref{bianchi-RS}
to remain valid, we need to assume that the hypersurfaces $P_t$ have
constant mean curvature and that the scalar curvature of $\bar{g}$ is constant
along each $P_t$. Note that the first variational formula
for volume of a hypersurface family (cf \cite{L}, Theorem 4) implies that the
formula $\dot{v} = \tr(L_t) v$ holds when $P_t$ are compact (oriented) and of
constant mean curvature. Finally, we also need to assume that
the divergence $\delta^{{\nabla}^t} \omega_t$ is constant on each $P_t$.
Under these conditions there is an analogous version of
Corollary \ref{cor-bianchi}. \end{rmk}

In the case of a gradient Ricci soliton with a $G$-equivariant potential,
we can make stronger statements as then $\omega_t =0$.

\begin{prop} \label{bianchi-GRS} Let $(\overline{M}, \bar{g})$ be a connected
Riemannian manifold with $\bar{g}$ of class $C^2$ and on which a compact Lie
group $G$ acts isometrically with cohomogeneity one. Let $u$ be a $C^3$
$G$-equivariant function on $\overline{M}$. Assume that $\overline{M}$ has
a special orbit with dimension smaller than that of a principal orbit.
If Eq.$($\ref{GRS-XY}$)$ holds on $M_0$, then so does $($\ref{GRS-XN}$)$.

If in addition $u$ satisfies $($\ref{altGRS-NN}$)$ on $M_0$, then $($\ref{GRS-NN}$)$
automatically holds.
\end{prop}

\noindent{\bf Proof.} We make use of the computations in the proof of
Proposition \ref{bianchi-RS}. For the first claim, we use (\ref{ricXN})
in (\ref{bianchi-1}) to get
$$0= \frac{\partial}{\partial t} \left(v( \overline{\ric}(X, N)) \right)
        - \dot{u} \, v  \overline{\ric}(X, N).$$
The set of $t \in {\rm int}(I)$ at which  $\overline{\ric}(X, N)$ does not vanish
is a disjoint union of open intervals. Over each such interval, one has
$| v \overline{\ric}(N, X)| = C \exp(u)$ for some positive constant $C$.
But $u(t)$ is everywhere defined and finite. So there is only one interval
and its endpoints must correspond to special orbits. If one of these has
dimension smaller than that of a principal orbit, we would have a contradiction.
Hence we conclude that  $\overline{\ric}(X, N)$ vanishes everywhere, which is
equivalent to the first claim.

For the second claim, we combine (\ref{altGRS-NN}) with (\ref{bianchi-2}) to
get
$$ N\left(v^2(\overline{\ric}(N, N) + \ddot{u} + \frac{\epsilon}{2})\right)
      =2v^2 \left(\overline{\ric}(N, N) + \ddot{u} + \frac{\epsilon}{2}\right) \dot{u}.$$
We can then argue exactly as for the first claim. \ \ \ \ \qed

\begin{rmk}
By Remark \ref{gen2} there are clearly analogous conclusions for the setting
of Remark \ref{gen} where we assume that there is a potential which is constant
on the fibres of the function $t$.
\end{rmk}

\section{\bf A Class of Equidistant Hypersurface Families}

In this section we consider the gradient K\"ahler Ricci soliton analogue
of the Einstein equations for the hypersurface families studied in \cite{WW},
\cite{DW}, \cite{Wa}, and \cite{CGLP}. Of course the hypersurfaces here need
not be homogeneous; in fact they provide an illustration of Remark \ref{gen}.

To fix notation, let $(V_i^{2n_i}, J_i, h_i), 1 \leq i \leq r,$ be respectively
compact K\"ahler-Einstein Fano manifolds with real dimension $2n_i$ and first Chern
class $c_1(V_i, J_i)= p_i a_i$, where $p_i$ are positive integers,
$a_i \in H^2(V_i; \Z)$ are indivisible classes, and the K\"ahler metric
$h_i$ is normalized by the condition $\ric(h_i)= p_i h_i.$ For $q=(q_1, \cdots, q_r)$
with $q_i \in \Z \setminus \{0\}$, let $P_q$ denote the principal $U(1)$-bundle over
$V:=V_1 \times \cdots \times V_r$ with Euler class $\sum_{1}^{r} q_i \pi_i^* a_i$,
where $\pi_i : V_1 \times \cdots \times V_r \rightarrow V_i$ is the projection
onto the $i$th factor. Denote by $M_0$ the product ${\rm int}(I) \times P_q$ for some
interval $I$.

Next let $\theta$ be the principal $U(1)$ connection on $P_q$ whose curvature is
$\Omega:= \sum_{1}^{r} q_i \pi_i^* \eta_i$ where $\eta_i$ is the K\"ahler
form of the metric $h_i.$ Using it we introduce the $1$-parameter family of
metrics $g_t$ on $P_q$ given by
\begin{equation} \label{metricfamily}
       g_t := f(t)^2 \theta \otimes \theta + \sum_{i=1}^{r} \, g_i(t)^2 \pi_i^*h_i
\end{equation}
where $f$ and $g_i$ are sufficiently smooth functions on $I.$
Note that each $g_t$ makes the bundle projection map into a Riemannian submersion
with totally geodesic fibres, and the curvature form $\Omega$ is parallel with
respect to any of the product metrics on the base.

We observe that $U(1)$ acts via isometries of $g_t$ on the right of $P_{q}$
for all $t$. This is important for two reasons. First it gives rise to a moment
map which allows us to introduce a reparametrization that simplifies the
Ricci soliton equations. Second if we choose the $1$-form $\overline{\omega}$
to be of the form
$$  \xi(t) dt +  \omega_t, \ \  \mbox{\rm where}   \ \ \ \omega_t = \lambda(t) \theta, $$
then it follows that $\delta^* \omega_t = \frac{1}{2} {\mathscr L}_{\lambda(t) f(t)^{-2} Z} g_t =0$
where $-Z$ is the Killing field generated by the right $U(1)$ action.
Notice that $\theta(Z)=1$ and $g_t(Z, Z) = f(t)^2.$
(The above form for $\overline{\omega}$ is a natural choice because when
all the $V_i$ are in addition homogeneous, then $P_q$ is homogeneous,
and for most choices of $q$, the trivial summand $\p_0$ in the isotropy
representation of $P_q$ is one-dimensional (spanned by $Z$)).

The metric $\bar{g}= dt^2 + g_t$
on ${\rm int}(I) \times P_q$ is easily seen to be hermitian with respect to
the complex structure $\overline{J}$ obtained by lifting the product complex
structure of the base to the horizontal spaces of $\theta$ and letting
$\overline{J}(N)= -f(t)^{-1}Z.$

We can now write down the Ricci soliton equations as in \cite{WW}. Now
$\tr(L_t) = \frac{\dot{f}}{f} + \sum_{i=1}^r 2n_i \frac{\dot{g_i}}{g_i}, $
is constant on $\{t\}\times P_q$. Also, using the fact that $g_t$ is a
Riemannian submersion with totally geodesic fibres, one sees that
$\delta^{\nabla^t} L_t = 0.$ Therefore we have $\overline{\ric}(N, X) = 0$
for all $X \in TP_q$. By the argument in Remark \ref{commondecomp}, we
are reduced to the case of a gradient Ricci soliton, i.e., the Ricci
soliton equation implies that $\lambda(t) \theta$, viewed as a $1$-form
on $\overline{M}$, is a Killing field and hence can be subtracted
off from $\overline{\omega}$.

The gradient Ricci soliton equation on ${\rm int}(I) \times P_q$ is then
the following system of equations (writing $\dot{u}$ for $\xi$ as in \S 1)

\begin{equation} \label{NN}
\frac{\ddot{f}}{f} + \sum_{i=1}^{r} \, 2n_{i}\frac{{\ddot{g}}_i}{g_i}
         - \ddot{u} = \frac{\epsilon}{2},
\end{equation}

\begin{equation} \label{UU}
\frac{\ddot{f}}{f} + \sum_{i=1}^{r} \, 2n_i \frac{\dot{f} \dot{g}_i}{fg_i}
- \frac{\dot{u}\dot{f}}{f}
-\sum_{i=1}^{r} \frac{n_i q_i^2}{2} \frac{f^2}{g_i^4}= \frac{\epsilon}{2},
\end{equation}

\begin{equation} \label{XX}
\frac{\ddot{g}_i}{g_i} - \left(\frac{\dot{g}_i}{g_i}\right)^2 + \frac{\dot{f} \dot{g}_i}{f g_i} +
\sum_{j=1}^r 2n_j \frac{\dot{g}_i \dot{g}_j}{g_i g_j}- \frac{\dot{u}\dot{g}_i}{g_i}
-\frac{p_i}{g_{i}^{2}} +  \frac{q_i^2 f^{2}}{2g_{i}^{4}}= \frac{\epsilon}{2},
\ \ \ \ 1 \leq i \leq r.
\end{equation}

As in \cite{WW} we introduce the moment map coordinate $s$
defined by $ds = f(t)dt$ and let
\begin{equation} \label{newvar}
\alpha(s):=f(t)^2, \ \ \ \ \beta_i(s):=g_i(t)^2, \ \ \ \ \varphi(s):= u(t).
\end{equation}
We shall denote by ${}^{\prime}$ differentiation with respect to $s$.
In addition, let us set
\begin{equation} \label{vol}
 v:= \prod_{i=1}^r g_i^{2n_i} = \prod_{i=1}^r \beta_i^{n_i}.
\end{equation}
Then the above system becomes
\begin{equation} \label{t-eqn}
\frac{1}{2} \alpha^{\prime\prime}+ \frac{1}{2}\alpha^{\prime}(\log v)^{\prime}
+ \alpha \sum_{i=1}^{r} \, n_{i}\left(\frac{\beta_{i}^{\prime\prime}}{\beta_i}
- \frac{1}{2} \left(\frac{\beta_i^{\prime}}{\beta_i}\right)^2 \right)
-\alpha \varphi^{\prime \prime} - \frac{1}{2} \alpha^{\prime} \varphi^{\prime} =\frac{\epsilon}{2},
\end{equation}

\begin{equation} \label{fibre-eqn}
\frac{1}{2}\alpha^{\prime\prime}+ \frac{1}{2}\alpha^{\prime}(\log v)^{\prime}
 -\frac{1}{2} \alpha^{\prime} \varphi^{\prime}-
 \frac{\alpha}{2} \sum_{i=1}^{r} \frac{ n_i q_i^2}{\beta_i^2} = \frac{\epsilon}{2},
\end{equation}

\begin{equation} \label{base-eqn}
\frac{\alpha^{\prime}}{2} \frac{ \beta_i^{\prime}}{\beta_i} +
\frac{\alpha}{2} \left(\frac{\beta_i^{\prime \prime}}{\beta_i} -
\left(\frac{\beta_i^{\prime}}{\beta_i} \right)^2 \right)
+ \frac{\alpha}{2}  \frac{\beta_i^{\prime}}{\beta_i} (\log v)^{\prime}
-\frac{\alpha}{2}  \frac{\beta_i^{\prime}}{\beta_i} \varphi^{\prime}
- \frac{p_i}{\beta_i} + \frac{q_i^2\alpha}{2\beta_i^2}= \frac{\epsilon}{2}.
\end{equation}

The first integral (\ref{altGRS-NN}) can be written as
\begin{equation} \label{integral}
\alpha \varphi^{\prime\prime} + \alpha^{\prime} \varphi^{\prime}
  + \alpha \varphi^{\prime}(\log v)^{\prime} - \alpha (\varphi^{\prime})^2 -
  \epsilon \varphi = c
\end{equation}
for some constant $c$. Equating (\ref{t-eqn}) with (\ref{fibre-eqn}) we
obtain
\begin{equation} \label{Jinvar}
 \varphi^{\prime\prime} = \sum_{i=1}^r \, n_i \left(\frac{\beta_i^{\prime \prime}}{\beta_i}
     -\frac{1}{2}\left(\frac{\beta_i^{\prime}}{\beta_i} \right)^2
     + \frac{1}{2}\frac{q_i^2}{\beta_i^2}\right).
\end{equation}

In view of the results in \cite{WW} and \cite{Wa} (Theorems 3.1 and 3.2)
we may look for solutions where
\begin{equation} \label{mu-i}
   \mu_i := \frac{\beta_i^{\prime \prime}}{\beta_i}
     -\frac{1}{2}\left(\frac{\beta_i^{\prime}}{\beta_i} \right)^2
     + \frac{1}{2}\frac{q_i^2}{\beta_i^2} =0, \ \ \  \ 1 \leq i \leq r.
\end{equation}
This forces $\varphi$ to be a linear function in $s$. In the following we
will show that this leads to gradient K\"ahler Ricci solitons generalizing
those constructed in \cite{Ko}, \cite{Ca1}, \cite{Ca2}, and \cite{FIK}.
As this analysis parallels that in \cite{WW}, we shall be brief and only
emphasize the necessary additional considerations.

We recall that setting all $\mu_i$ to be $0$ is equivalent to the
geometric condition that the Riemann curvature tensor of $\bar{g}$ is
fully invariant under the action of the complex structure $\bar{J}$
(cf Corollary 7.5 in \cite{WW}). There are two types of solutions to
$\mu_i=0:$ either $\beta_i$ is a quadratic polynomial in $s$ of the
form $A_0(s+s_0)^2 - \frac{q_{i}^2}{4A_0}$ with $A_0 \neq 0$ or it is
a linear polynomial of the form $\pm q_i(s+ \sigma_i)$. Furthermore,
in the latter case, the choice of the minus sign for all $i$
corresponds to the K\"ahler condition for the metric $\bar{g}$ with
respect to the complex structure $\bar{J}$. We recall from \cite{G} or
\cite{FIK} that for gradient Ricci solitons where the metric is
K\"ahler, the vector field $X = {\rm grad}\, u$ is automatically an
infinitesimal automorphism of $\bar{J}$.

Accordingly, let us set
\begin{equation} \label{betaform}
 \beta_i :=-q_i(s + \sigma_i), \ \ \ \ \ \ {\rm and} \ \ \ \ \varphi:=\kappa_1(s+ \kappa_0),
\end{equation}
where $\sigma_i$ and $\kappa_i$ are real constants to be determined.
Substituting these into (\ref{base-eqn}) leads to
\begin{equation} \label{alpha-KE}
   \alpha^{\prime} + \alpha((\log v)^{\prime} -\kappa_1)  = \epsilon s + E^*
\end{equation}
and the consistency condition
\begin{equation} \label{consist1-KE}
E^* = \epsilon \sigma_i -\frac{2p_i}{q_i}, \ \ \ \ \ \  1 \leq i \leq r.
\end{equation}
One can now directly check that (\ref{t-eqn}) and hence (\ref{fibre-eqn}) automatically
hold, while (\ref{integral}) is consistent with (\ref{base-eqn}) provided that
\begin{equation} \label{consist2-KE}
 E^* = \epsilon \kappa_0 +\frac{c}{\kappa_1}.
\end{equation}
Integrating (\ref{alpha-KE}) gives
\begin{equation} \label{alphaform1}
 \alpha = v^{-1} e^{\kappa_1 s} \int \, (\epsilon s + E^*)\,  e^{-\kappa_1 s} v \; ds.
\end{equation}
Combined with the earlier expression (\ref{betaform}) we have a family of
{\em explicit} solutions to the equations.

In order to construct complete examples from the above local solutions, we need to
analyse the smoothness conditions when we compactify an end of our cohomogeneity
one manifold $M_0 = {\rm int}(I) \times P_q$ by adding a compact submanifold at $t=0$.
We may assume without loss of generality that $s=0$ when $t=0$.
There are then two possibilities. We can add $V_1 \times \cdots \times V_r$, which
corresponds to letting $\alpha$ go to zero as $s$ tends to $0$.
If $V_1,$ say, is a projective space $\C \PP^{n_1}$
then we can add $V_2 \times \cdots \times V_r$ ($r \geq 2$), which corresponds to
letting both $\alpha$ and $\beta_1$ go to zero as $s$ tends to $0$. Notice that
the former case may be regarded as a special case of the latter case by
allowing $n_1$ to be zero.

When we compactify $M_0$ by adding $V_1 \times \cdots \times V_r$ at $t=0$,
the smoothness conditions for the metric $\bar{g}$ can be deduced by
the methods in section 1 of \cite{EW}. One concludes that $f$ should be smooth
and odd in $t$ with $\dot{f}(0)=1$ and that $g_i$ should be smooth and even
with $g_i(0) \neq 0$. (We have used here the condition $\theta(Z)=1$.)
However, by Lemma \ref{DTKz}, we actually only need to check these conditions
up to order $2$ and that the functions have a finite third derivative at $t=0$.
>From $f(t)^2 = \alpha(s)$ it follows that if $\alpha$ is thrice differentiable
with $\alpha(0)=0$ and $\alpha^{\prime}(0)=2,$ then $f$ is odd up to order
$2$ and thrice differentiable with $\dot{f}(0)=1.$ Using these properties
in conjunction with the differentiability of  $\varphi$ and $\beta_i$ and
(\ref{newvar}), it follows easily that $u$ and $g_i$ are thrice differentiable
and even (up to order $2$). Also, $\beta_i(0) > 0$ gives $g_i(0) \neq 0.$

If $V_1 = \C \PP^{n_1}$ and we compactify $M_0$ by adding the submanifold
$V_2 \times \cdots \times V_r$, then the above analysis remains valid
except that we must suitably change the conditions on $\beta_1, g_1(t)$
in order that $\bar{g}$ is smooth. Here the key point is that the distance
spheres in the normal bundle of $V_2 \times \cdots \times V_r$ must become
round as $t$ tends to $0$. Recall that in the Hopf fibration
$S^{2n_1 +1} \rightarrow \C \PP^{n_1}$ where the sphere has constant
curvature $1$, the submersed metric has Einstein constant $2(n_1 +1)$.
Since we have chosen $\dot{f}(0)=1$ and $\ric(h_1) = (n_1 +1) h_1$, this
implies that $\dot{g_1}(0)^2 = 1/2.$ The relation $2g_1 \dot{g_1} = \beta_1^{\prime} f$
then forces us to choose $\beta_1^{\prime}(0) = 1$, i.e., $q_1 = -1$.
In other words, if we let $q_1 =-1$ and $\beta_1(s) = s$, then $g_1(t)$ will
be smooth and odd (up to order $3$) with $\dot{g_1}(0) = \pm 1/\sqrt{2},$
which will ensure that $\bar{g}$ is smooth.

\medskip
\noindent{\bf \em Steady Solitons}

\medskip
We let $\epsilon =0$ and $V_1 = \C \PP^{n_1}$ with $n_1 \geq 0$.
Let $\beta_1(s) =s,$ i.e., we take $q_1 =-1$ and $\sigma_1=0$.
 In order to have $\beta_i(s) > 0$ on $[0, +\infty)$,
$2 \leq i \leq r,$ we need to assume $-q_i > 0$ and $\sigma_i > 0$.
The consistency conditions (\ref{consist1-KE}) and (\ref{consist2-KE}) lead to
$$ 2n_1 + 2 =E^* = \frac{c}{\kappa_1}
=-\frac{2p_i}{q_i},    \ \ \ 2 \leq i \leq r, $$
except in the case $r=2, n_2 =0$ when there is no third equality.
The manifold $\overline{M}$ is then $\C^{n_1 +1}.$
In any case it follows that
\begin{equation} \label{alphaform2}
 \alpha(s) = \frac{(2n_1+2)e^{\kappa_1 s}}{\prod_{i=1}^r (s+\sigma_i)^{n_i}}
         \int_0^s e^{-\kappa_1 x} \prod_{i=1}^r (x + \sigma_i)^{n_i} dx.
\end{equation}
Notice that $v = \mu s^{n_1} +$ higher powers of $s$, where $\mu$ is a
constant. It follows from (\ref{alphaform1}) that $\alpha(s) =
\frac{E^*}{n_1+1}s + \ldots$, so the above constraint $E^* = 2n_1 + 2$
guarantees that $\alpha^{\prime}(0)=2$, as required by the collapsing.
It is also clear that $\alpha(0)=0$ and, from (\ref{alphaform2}), that
$\alpha(s) > 0$ when $s>0$ since we've chosen $\sigma_i$ positive for
$i \geq 2$. This, together with our choice of $q_1$ above, shows that
 the metric extends smoothly to the compactified end at $s=0$.

To ensure that we get a complete metric $\bar{g}$, recall that the geodesic
distance
$$ t = \int_0^s \frac{dx}{\sqrt{\alpha(x)}},$$
and so we need to ensure that the integral diverges as $s \rightarrow +\infty$.
If $\kappa_1> 0,$ then the integral in (\ref{alphaform2}) is bounded and so
$\alpha(s)$ grows exponentially. The geodesic distance would then be bounded.
On the other hand, if $\kappa_1 < 0$, the integral in (\ref{alphaform2})
grows like $\frac{1}{-\kappa_1} s^{n_1 + \cdots + n_r} e^{-\kappa_1 s}$.
So $\alpha(s)$ is asymptotic to a positive constant, and hence the geodesic
distance is unbounded and the metric is complete.
 The $\kappa_1 = 0$ case is the Einstein case, and
$\alpha(s)$ grows like $s$. So again the geodesic distance is unbounded.

For $\kappa_1 < 0$,
we see that (suppressing multiplicative constants)
for large $t$, $f(t)$ is $O(1)$,  $g_i(t) \sim t^{1/2}$,
$u(t) \sim t,$ and the volume of the hypersurfaces grow like  $t^{(n-1)/2}$.
(Recall that $n$ is the dimension of $P_q$.) In other words, the hypersurfaces
$\{t\} \times P_q$ are asymptotically circle bundles whose fibres have
approximately constant circumference. This kind of behaviour is often referred to
as cigar-paraboloid asymptotics (inspired by the Hamilton-Witten cigar soliton
in complex dimension one).
One can calculate in our examples that the sectional curvatures
of $\bar{g}$ decay at least as fast as $t^{-1}$ as $t$ becomes large.

The free parameters in this construction are $\sigma_i > 0$ and $\kappa_1 \leq 0$.
(The constant $\kappa_0$ represents an inherent ambiguity of the potential function
$u(t)$.) If we multiply the metric $\bar{g}$ by a positive constant,
then $\epsilon$ becomes divided by the constant, as the Hessian of $u$
and the Ricci tensor are unchanged. Therefore, we obtain an $r$ (resp. $r-1$)
parameter family of K\"ahler Ricci soliton solutions if we compactify $M_0$ at one end
by adding $V_1 \times \cdots \times V_r$ (resp. $V_2 \times \cdots \times V_r$).

Note finally that if we let $\kappa_1=\ddot{u}(0)$ tend to $0$ we obtain
the convergence of non-trivial K\"ahler Ricci solitons to Ricci-flat K\"ahler metrics.

\medskip

\noindent{\bf \em Expanding Solitons}

\medskip

We assume $\epsilon =1$ to factor out homothety. As in the case of steady solitons,
we need $-q_i > 0, \sigma_i > 0$ for $2 \leq i \leq r$ and $q_1 = -1$. The
consistency conditions (\ref{consist1-KE}) and (\ref{consist2-KE}) become
$$  \kappa_0 + \frac{c}{\kappa_1} = E^* = 2n_1 +2 = -\frac{2p_i}{q_i} + \sigma_i,
\ \ \ 2 \leq i \leq r.$$
Hence $\sigma_i$ are no longer free parameters and we need to have
$$  -q_i (n_1 + 1) > p_i,$$
except when $r=2, n_2 =0.$ We now have
\begin{equation}
 \alpha(s) = \frac{e^{\kappa_1 s}}{\prod_{i=1}^r (s+\sigma_i)^{n_i}}
         \int_0^s (x + 2n_1 +2) e^{-\kappa_1 x} \prod_{i=1}^r (x + \sigma_i)^{n_i} dx.
\end{equation}
As before, we need $\kappa_1 \leq 0$ for the completeness of $\bar{g}$.
Indeed when $\kappa_1 < 0$, $\alpha(s)$ grows like $-\frac{1}{\kappa_1} s$
and when $\kappa_1 = 0$, $\alpha(s)$ grows like $s^2$, so the geodesic
distance $t$ is unbounded. In the
case $\kappa_1 < 0$, it follows that
(suppressing multiplicative constants) for large $t$, $f(t) \sim t$,
$g_i(t) \sim t,$ $u(t) \sim t^2$ and the hypersurfaces $\{t\} \times P_q$
have Euclidean volume growth, i.e.,  each  $\overline{M}$ has an asymptotically
conical end. Hence we obtain a one-parameter family of soliton
solutions  which converge to a complete K\"ahler-Einstein metric with negative
scalar curvature. As before, our choice of $E^*$ and $q_1$ guarantees smooth extension
over the compactified end at $s=0$. One may check that the sectional
curvatures of the metric decay like $t^{-2}$ as $t$ becomes large.

We may summarise the above discussion as

\begin{thm} \label{completeGKRS}
Let $(V_i, J_i, h_i),  1 \leq i \leq r, \ r\geq 2,$ be Fano K\"ahler-Einstein
manifolds with complex dimension $n_i$ and first Chern class $p_i a_i$
where $p_i > 0$ and $a_i$ are indivisible classes in $H^2(V_i, \Z)$. Let $V_1$
be $\C \PP^{n_1}$ with normalised Fubini-Study metric and assume that $n_1 \geq 0$.
Let $P_q$ denote the principal $S^1$ bundle over $V_1 \times \cdots \times V_r$
with Euler class $-\pi_1^*(a_1) + \sum_{i=2}^r q_i \, \pi_i^*(a_i)$.
\begin{enumerate}
\item[$($i$)$] If $-q_i(n_1 +1) = p_i$ for all $2 \leq i \leq r,$ then there
is an $(r-1)$-parameter family of non-trivial complete steady gradient K\"ahler-Ricci
solitons on the underlying space of the corresponding complex $\C^{n_1 +1}$
vector bundle over $V_2 \times \cdots \times V_r$.
\item[$($ii$)$] If $-q_i(n_1 +1) > p_i,$ for all $ 2 \leq i \leq r$, then there is a
$1$-parameter family of non-trivial complete expanding gradient K\"ahler-Ricci
solitons on the corresponding complex $\C^{n_1 +1}$ vector bundle over
$V_2 \times \cdots \times V_r$.
\end{enumerate}
In both cases, the K\"ahler metric on the bundle has a circle of isometries
and the soliton potential $u$ can be chosen to be constant on the distance
sphere-subbundles. If we  let the value of $\ddot{u}$ at the zero section
tend to $0$ and fix the rest of the  parameters, then the soliton metrics
converge to a complete K\"ahler-Einstein metric. \, \,  \qed
\end{thm}

\medskip

\begin{rmk} \label{history1}
In (i) of the above theorem, the case $n_1=0, r=2$ with $V_2 = \C \PP^n$
was obtained in \cite{Ca1} and \cite{ChV} (cf Proposition 5), as was the
case $n_2 = 0, r=2$ with $V_1 = \C \PP^{n}$ where the manifold is $\C^{n+1}.$
Cao further made the important observation that this $1$-parameter family
of examples contain ones with positive sectional curvature.

In (ii) of the above theorem, the case $n_2 =0, r=2$ with $V_1 = \C \PP^{n}$
was  obtained in \cite{ChV} (cf Proposition 3) and \cite{Ca2}.
Cao again showed that this example contained solutions with positive
sectional curvature. The case $n_1 = 0, r=2$ with $V_2 = \C \PP^{n}$ was
obtained independently in \cite{ChV} (cf Proposition 5) and \cite{FIK}.
That $V_2$ can be replaced by any Fano K\"ahler-Einstein manifold was noted
in \cite{PTV} (cf Theorem 2).
\end{rmk}

\begin{rmk}
The above arguments also work in the case of {\em expanding} solitons
if some of the $p_i$ are allowed to be nonpositive (i.e., the K\"ahler-Einstein
manifolds have $c_1 \leq 0$). When the base consists of a single K\"ahler-Einstein
factor this was observed in \cite{PTV} (cf Theorem 1).

More precisely, for a base factor that has negative first Chern class,
we need to assume that $c_1(V_i, J_i) = p_i a_i$ where $p_i$ is a negative integer
and $a_i$ is an indivisible class in $H^2(V_i; \Z)$, and for a Calabi-Yau
base factor, we assume that the K\"ahler form $\eta_i$ of the metric $h_i$
is $2\pi$ times an {\em integral} cohomology class $a_i$. Then the conditions on
the Euler class of the $U(1)$ bundle $P_q$ are that $q_1=-1$ \, ($V_1 = \C \PP^{n_1},
n_1 \geq 0$ is still the collapsing factor), $q_i < 0$ for the non-positive KE factors,
and $-q_i(n_1 +1) > p_i$ for the Fano KE factors. So again we obtain a $1$-parameter
family of K\"ahler Ricci solitons.
\end{rmk}

\noindent{\bf \em Shrinking Solitons}

\medskip

We again set $\epsilon = -1$ to factor out homothety. We first assume $r \geq 2$
and let $V_1 = \C \PP^{n_1}$ ($n_1 \geq 0$) with the normalised Fubini-Study
metric. Then the consistency conditions (\ref{consist1-KE}) and (\ref{consist2-KE})
become
\begin{equation} \label{consist3}
 2n_1 + 2 = E^*= \frac{c}{\kappa_1}-\kappa_0 = -\sigma_i - \frac{2p_i}{q_i},
 \, \, \, 2 \leq i \leq r.
\end{equation}
Observe that if we ensure that $v>0$ in (\ref{alphaform1}), except
at endpoints of $I$, the factor $E^* -s = 2n_1 +2 -s$ in the integrand
could still make $\alpha$ become non-positive. We therefore have to consider
two separate cases, corresponding to whether or not this happens at a finite
endpoint.

First let us consider solutions defined on a finite interval $[0, s_*],$ where at
$s=s_*$ we need to put in a compactifying submanifold as well.  In this case,
let $r \geq 3$ and set $V_r = \C \PP^{n_r}$ ($n_r \geq 0$), equipped with the
normalised Fubini-Study metric. The compactifying submanifold at $s=s_*$ becomes
$V_1 \times \cdots \times V_{r-1}$.

If we now apply the smoothness conditions at $s=s_*$, we obtain $q_r = 1$
and $\beta_r = s_*-s$. The consistency condition (\ref{consist1-KE}) implies
that $s_*=2(n_1 + n_r +2)$. In order to have $\beta_i > 0$ on $[0, s_*]$,
the inequalities
$$ -(n_1 +1) q_i < p_i, \, \, \, \, (n_r +1)q_i < p_i, \, \, \, \, 2 \leq i \leq r-1 $$
must hold. These relations mean that $\sigma_i, s_* + \sigma_i$ both have the opposite
sign to $q_i$ for $2 \leq i \leq r-1$.
Now $\alpha$ becomes
$$\alpha(s) = \frac{e^{\kappa_1 s}}{\prod_{i=1}^r (s-2n_1-2-\frac{2p_i}{q_i})^{n_i}}
         \int_0^s (2n_1 +2-x  ) e^{-\kappa_1 x}
         \prod_{i=1}^r \left(x -2n_1 -2 -\frac{2p_i}{q_i} \right)^{n_i} dx, $$
where we have cancelled a factor of $\prod_i (-q_i)^{n_i}$ from the numerator
and denominator. By examining this formula, we see that in order for
$\alpha(s_*)= 0$ and $\alpha(s) > 0$ on $(0, s_*)$ it is necessary and
sufficient that the integral
\begin{equation} \label{Futaki}
 \mathcal I :=\int_{-n_1 -1}^{n_r +1} \, e^{-2\kappa_1(x+n_1 +1)}
      \prod_{i=1}^{r} \, \left(x-\frac{p_i}{q_i}\right)^{n_i}\, x\, dx =0.
\end{equation}
We have therefore deduced

\begin{prop} \label{koiso}
Let $(V_i, J_i, h_i),  1 \leq i \leq r, \ r\geq 3,$ be Fano K\"ahler-Einstein
manifolds with complex dimension $n_i$ and first Chern class $p_i a_i$
where $p_i > 0$ and $a_i$ are indivisible classes in $H^2(V_i; \Z)$. Let $V_1$
and $V_r$ be complex projective spaces with normalised Fubini-Study metric.
Let $P_q$ denote the principal $S^1$ bundle over $V_1 \times \cdots \times V_r$
with Euler class
$-\pi_1^*(a_1) + \sum_{i=2}^{r-1} q_i \, \pi_i^*(a_i) +\pi_r^*(a_r)$.

Suppose in addition that $-(n_1 +1)q_i < p_i$ and $ (n_r+1)q_i < p_i$
for all $2 \leq i \leq r-1$. Then there is a compact shrinking
gradient K\"ahler Ricci soliton structure on the space $\overline{M}$ obtained from
$P_q \times_{S^1} \C\PP^1$ by blowing $P_q$ down to $V_2 \times \cdots \times V_r$
at one end and to $V_1 \times \cdots \times V_{r-1}$ at the other end iff for some
$\kappa_1 \in \R$, the integral in $($\ref{Futaki}$)$ vanishes. The Ricci soliton
is K\"ahler-Einstein if $\kappa_1=0$ and is otherwise non-trivial. \ \ \ \ \qed
\end{prop}

\begin{rmk} \label{history2}
(i) The examples of Koiso \cite{Ko}, Cao \cite{Ca1}, and Chave-Valent \cite{ChV}
of $\C \PP^1$-bundles over complex projective space correspond to the situation
when $r=3,$ $V_2 = \C \PP^{n_2}$, and $n_1 = n_3 = 0$. Here $q_2$ must satisfy
$0 < |q_2| < n_2 +1$. Other similar examples with orbifold singularities were
constructed in \cite{FIK}. Actually, it was already observed earlier in \cite{G}
that $V_2$ can be any Fano K\"ahler-Einstein manifold.

(ii) The integral $\mathcal I$ is a special case of the new holomorphic invariant 
introduced in \cite{TZ2}. When $\kappa_1=0$, it becomes
  the Futaki character
for the first Chern class of $(\overline{M}, \overline{J})$ evaluated on the
real holomorphic vector field $\overline{J}(Z)=fN$.

(iii) Note also that if $n_1 =n_r$ we have (for fixed $p_i$ and $n_i$) the relation
$${\mathcal I}(-\kappa_1, -q) =
   (-1)^{1+\sum_i n_i} e^{4\kappa_1(n_1+1)} {\mathcal I}(\kappa_1, q).$$
In particular, when $n_1=n_r$, if ${\mathcal I}$ vanishes at $\kappa_1 =\varrho$
for the K\"ahler manifold determined by $q=(q_1, \cdots, q_r),$  then
$\mathcal I$ vanishes at $\kappa_1=-\varrho$ for the K\"ahler manifold determined
by $-q$ and interchanging $V_1$ and $V_r$ . (These manifolds are related by a
diffeomorphism which reverses orientation along the fibres.)
\end{rmk}

We will now examine the asymptotics of the integral in (\ref{Futaki})
as $|\kappa_1|$ becomes large. We shall use
the identity
\begin{equation}\label{expidentity}
  \int_0^s e^{-\kappa_1 x} x^m dx = \frac{m!}{\kappa_1^{m+1}}
       \left(1 - e^{-\kappa_1 s} \sum_{j=0}^m \frac{(\kappa_1 s)^j}{j!}
\right).
\end{equation}
In particular, note that if ${\sf P}(x)$ is a polynomial
then for $\kappa_1 >> 0$
\[
\int_{0}^{s} e^{-\kappa_1 x} {\sf P}(x) \;dx \sim \frac{b_m m!}{\kappa_1^{m+1}}
\]
where $b_m$ is the lowest nonzero coefficient of $\sf P$.
Letting $x = s - \tilde{x}$, we see that for $\kappa_1 << 0$
\[
\int_{0}^{s} e^{-\kappa_1 x} {\sf P}(x) \; dx \sim e^{-\kappa_1 s}
\frac{b_{\tilde{m}} \tilde{m}! }{(- \kappa_1)^{\tilde{m}+1}}
\]
where ${\sf P}(x) = b_{\tilde{m}} (s -x)^{\tilde{m}} + \; $ higher powers
of $(s-x)$. Hence the asymptotic signs of the integral
are given by the signs of $b_m$ and $b_{\tilde{m}}$ respectively.

\medskip
If we now substitute $y = 2(x + n_1 +1)$ in  the integral in (\ref{Futaki}),
we obtain
\[
{\mathcal I}=2^{-(n_1 +\ldots +n_r +2)} \int_{0}^{s_*} e^{-\kappa_1 y} (y - 2n_1 -2)
\prod_{i=1}^{r} \left( y - 2n_1 -2 - \frac{2p_i}{q_i} \right)^{n_i} \; dy.
\]
Equivalently,
\begin{equation} \label{Futaki2}
{\mathcal I} = 2^{-(n_1 +\ldots +n_r +2)}
\int_{0}^{s_*} e^{-\kappa_1 y} y^{n_1} (y-s_*)^{n_r}
(y -2n_1 -2) \prod_{i=2}^{r-1}(y + \sigma_i)^{n_i} \;dy
\end{equation}
where we have used (\ref{consist3}),
$p_1 = n_1 +1, p_r = n_r+1, q_1 =-1, q_r =1,$ and $ s_*=2(n_1 +n_r +2)$.
Now the discussion following (\ref{expidentity}) shows that
\begin{equation} \label{asympplus}
{\mathcal I} \sim \frac{2^{-(n_1 + \ldots +n_r +1)} n_1 ! s_*^{n_r} (n_1 +1)}
{\kappa_{1}^{n_1+1}} \, (-1)^{n_r +1}\prod_{i=2}^{r-1} \sigma_i^{n_i} \;\;
{\rm for} \; \kappa_1 >> 0.
\end{equation}
In particular the asymptotic sign of $\mathcal I$ is that of $(-1)^{n_r +1}
\prod_{i=2}^{r-1} \sigma_{i}^{n_i}$.

Similarly, the asymptotic sign of $\mathcal I$ for $\kappa_1 << 0$ is that of
$(-1)^{n_r} \prod_{i=2}^{r-1} (s_* + \sigma_i)^{n_i}$. But as noted earlier,
$s_*+\sigma_i$ and $\sigma_i$ have the same sign for $2 \leq i \leq r-1$.
Hence the asymptotic signs are always opposite, and we obtain

\begin{thm} \label{koisoII}
All the compact K\"ahler manifolds $\overline{M}$
described in Theorem \ref{koiso} are Fano and
admit an explicit gradient K\"ahler Ricci soliton.    \ \ \ \ \ \ \  \qed
\end{thm}

We illustrate our discussion with some examples.

\begin{example}
Let us consider Ricci solitons on $\C\PP^1$-bundles over $\C \PP^2 \times \C \PP^2$.
We are therefore just collapsing a circle at each endpoint, rather than a
higher-dimensional sphere, so we take $r=4$ and $n_1=n_4=0$. Moreover
$n_2 = n_3=2$ and $p_2 =p_3=3$. We must choose $|q_i | < p_i$ for
$i=2, 3$. (This also implies that the first Chern class of the resulting
complex manifold is positive.) Let us take $(q_2, q_3)=(1, -2)$. Our integral
(\ref{Futaki}) becomes:
\begin{equation} \label{FutakiEx}
\int_{-1}^{1} e^{-2 \kappa_1(x+1)} (x-3)^2 (x + \frac{3}{2})^2 x \; dx.
\end{equation}
Using MAPLE, we find that if $\kappa_1=0$ then this integral equals 7.8 (so there
is no K\"ahler-Einstein metric).
If $\kappa_1 = \frac{1}{2}$, on the other hand, the integral is
approximately -0.7289, so there is a $\kappa_1 \in (0, \frac{1}{2})$ where
(\ref{FutakiEx}) vanishes and we have a shrinking Ricci soliton.

Note that more examples can be obtained by taking the base to be a
product of complex projective spaces. Since ${\rm SU}(n_i+1)$ acts
transitively on $\C\PP^n_i$ and  we have an additional isometric circle action
on the fibres, the resulting $S^2$ bundles are {\em toric}. There is, of course, a
general existence theorem for K\"ahler Ricci solitons for toric Fano varieties
due to Wang and Zhu \cite{WZ}. So the examples we get here are not new,
but the K\"ahler metrics are reasonably explicit.
\end{example}

For non-toric examples with inhomogeneous base
we may take some of the factors in the base to
be suitable Fermat hypersurfaces ${\mathcal F} (n,d)$, i.e., smooth
degree $d$ hypersurfaces in $\C \PP^{n+1}$. These have $c_1$
equal to $(n+2-d)$ times the generator of the second integral
cohomology group, so are Fano if $d < n+2$. They are known to admit
K\"ahler-Einstein metrics by the work of Siu \cite{Si} and Tian \cite{T1}
if $d=n, n+1,$ of Nadel \cite{Na} if $\frac{1}{2}(1+ n) \leq d \leq n+1$ ,
and of Tian \cite{T2} for the remaining cases.

\begin{example}
Let us take $r=4$ and $n_1 = n_4=0$, with $V_2$ and  $V_3$ equal
to the irrational Clemens-Griffiths three-fold ${\mathcal F}(3,3)$
\cite{ClGr}.
As $n_2 = n_3 = 3$, $c_1 ({\mathcal F}(3,3))$ is
twice the generator, and so $p_2=p_3=2$.

If we choose $(q_2, q_3)=(-1,-1)$ the integral (\ref{Futaki}) is
\[
\int_{-1}^{1} \; e^{-2 \kappa_1 (x+1)}(x+2)^6 x \; dx.
\]
If $\kappa_1 =0$, this is $\frac{1368}{7}$, so there is no
K\"ahler-Einstein metric. If $\kappa_1$ is large positive, we see from
(\ref{asympplus}) that the integral is $\sim -\frac{1}{2 \kappa_1}$,
so we deduce there is a positive value of $\kappa_1$ for which the
integral vanishes. We thus obtain a K\"ahler-Ricci soliton on a $\C
\PP^1$-bundle over a product of two copies of ${\mathcal F}(3,3)$.
\end{example}

\begin{example}
We next consider an example with blow-downs. We take $r=3$, $n_1=n_3=1$ and
$V_2$ to be ${\mathcal F}(4,3)$. Now $n_2 =4$ and
$c_1 ({\mathcal F}(4,3))$ is three times
the generator . We need $|2q_2| <p_2 =3$, so $q_2 = \pm 1$. If we choose
$q_2 =-1$ our integral becomes
\[
\int_{-2}^{2} e^{-2 \kappa_1 (x+2)} (x+2)(x + 3)^4 (x-2) x \; dx.
\]
For $\kappa_1=0$ this is $-\frac{7680}{7}$ while for $\kappa_1$ large positive
 it is $\sim \frac{2}{\kappa_1^2}$, so again we deduce the existence of a
soliton.
\end{example}

\medskip
Let us now consider noncompact complete shrinking solitons.
We look for solutions  defined on $[0, +\infty),$ so that
we only need to put in the compactifying submanifold $V_2 \times \cdots \times V_r$
at $s=0$. That this type of solution actually exists was first observed in \cite{FIK}.
In the following we describe the natural generalization of these examples.

As in the steady and expanding cases, we take $r \geq 2$ with $q_1 = -1, \sigma_1 =0$,
so $\beta_1(s)=s$. In order for $\beta_i(s) > 0$, $2 \leq i \leq r$,
we need to assume $-q_i > 0$ and $\sigma_i > 0$. The consistency
conditions (\ref{consist1-KE}) and (\ref{consist2-KE}) then lead to
the conditions
$$0 < -(n_1 +1)q_i < p_i, \,\,  2 \leq i \leq r.$$
The specific form for $\alpha(s)$ becomes
$$\alpha(s) = \frac{e^{\kappa_1 s}}{s^{n_1} \prod_{i=2}^r |q_i|^{n_i}(s+ \sigma_i)^{n_i}}
         \int_0^s  e^{-\kappa_1 x} \left((2n_1 +2-x) x^{n_1}
         \prod_{i=2}^r |q_i|^{n_i}(x + \sigma_i)^{n_i} \right) dx. $$
It is clear from this expression that unless $\kappa_1 > 0$, $\alpha(s)$ will
eventually become negative. Therefore, from now on we let $\kappa_1 > 0.$
Since $\sigma_i > 0$ by choice, the integrand is positive on $(0, 2n_1 +2)$
and negative if $x > 2n_1 +2.$ So the integral in $\alpha(s)$
is increasing on $(0, 2n_1 +2)$ and monotonically decreasing on
$(2(n_1 + 1), +\infty).$ In particular, $\alpha(s) > 0$ for all $s>0$ provided
we can show that the integral is asymptotically positive.

Observe as before that the term in braces in the integrand is a polynomial
$$ \Psi(x) = a_D x^D + a_{D-1} x^{D-1} + \cdots + a_{n_1}x^{n_1} $$
with $D = n_1 + \cdots + n_r +1$ and $a_D < 0, \, a_{n_1} > 0$. Using the
formula (\ref{expidentity})
we obtain
$$ s^{n_1} \prod_{i=2}^{r} (|q_i|(s+ \sigma_i))^{n_i} \alpha(s) =
   e^{\kappa_1 s} \int_0^s e^{-\kappa_1 x} \Psi(x) dx =
     e^{\kappa_1 s} \sum_{k=n_1}^D \,  \frac{k! a_k}{\kappa_1^{k+1}} -
     \sum_{k=n_1}^D \frac{k! a_k}{\kappa_1^{k+1}} \left(\sum_{j=0}^k \frac{(\kappa_1 s)^j}{j!}   \right).$$

Recall from our discussion of the steady and expanding cases
that if we are to have a complete metric $\alpha(s)$ cannot grow exponentially.
It follows from our last formula that
we must choose $\kappa_1 > 0$ so that
\begin{equation} \label{kappa1}
 \sum_{k=n_1}^D  \frac{k! \, a_k}{\kappa_1^{k-n_1}} = 0.
\end{equation}
This is certainly possible since $a_D < 0$ and $a_{n_1} > 0$.
In fact, such a $\kappa_1$ is unique. For the expression on the left of
(\ref{kappa1}) is $\chi (\frac{1}{\kappa_1})$ where the coefficients
of $\chi(x)$ are obtained from those of $\Psi(x)/x^{n_1}$ by
multiplying by positive constants. But by the definition of $\Psi$,
we see the roots of the polynomial $\Psi(x)/x^{n_1}$ are all
real, and moreover one root is positive and the rest are all negative.
Descartes's rule of signs now implies there is exactly one sign change in the
coefficients of $\Psi(x)/x^{n_1}$. Hence this is also true
for $\chi(x)$, and the rule of signs now implies $\chi$ has a unique
positive root.

With such a choice of $\kappa_1$, it follows that for large $s$
we have
$$ \alpha(s) \sim \frac{-a_D}{\kappa_1 \prod_{i=2}^r |q_i|^{n_i}} \, s
\sim \frac{s}{\kappa_1}$$
As  $\kappa_1$ is positive, we see that $\alpha(s)$ is positive for large
$s$, and hence as observed earlier, for all $s> 0.$ So our Ricci soliton
solution is defined on the whole interval $[0, \infty)$, as desired.
Furthermore, the geodesic distance $t \sim 2 \sqrt{\kappa_1} s^{1/2},$
so $f(t) \sim (2 \kappa_1)^{-1}t,$ and $g_i(t), u(t) \sim$ positive constants
times $t$ and $t^2$ respectively. Hence $\overline{M}$ has an asymptotically
conical end, just like the situation in the expanding case. In particular,
the metric is complete. As in the expanding case, the sectional curvatures
decay like $t^{-2}$. We have therefore deduced

\begin{thm} \label{FIK}
Let $(V_i, J_i, h_i),  1 \leq i \leq r, \ r\geq 2,$ be Fano K\"ahler-Einstein
manifolds with complex dimension $n_i$ and first Chern class $p_i a_i$
where $p_i > 0$ and $a_i$ are indivisible classes in $H^2(V_i, \Z)$. Let $V_1$
be $\C\PP^{n_1},  n_1 \geq 0,$ with normalised Fubini-Study metric and
let $P_q$ denote the principal $S^1$ bundle over $V_1 \times \cdots \times V_r$
with Euler class
$-\pi_1^*(a_1) + \sum_{i=2}^{r} \, q_i \, \pi_i^*(a_i) $.

Suppose in addition that $0< -(n_1 +1)\, q_i < p_i$ for all $2 \leq i \leq r$.
Then there is a complete shrinking gradient K\"ahler Ricci soliton structure on
the space $\overline{M}$ obtained from the line bundle $P_q \times_{S^1} \C$
by blowing the zero section down to $V_2 \times \cdots \times V_r.$
The Ricci soliton metric has an asymptotically conical end. \ \ \ \ \qed
\end{thm}

\begin{rmk} The examples in \cite{FIK} correspond to taking $r=2, n_1 =0$
and $V_2$ to be a complex projective space. As noted there, the case
$r=2, n_2=0$ corresponds to flat $\C^{n_1+1}$ as a shrinking soliton.
Of course the essentials of the above analysis are similar to those
in \cite{FIK}.
\end{rmk}

As discussed in \S 0, the Ricci flow of a soliton with vector field
$X$ is a combination of rescaling by $(1 + \epsilon \tau)$ and pulling
back by diffeomorphisms $\psi_{\tau}$, where $\psi_{\tau}$ integrate
the field $Y_{\tau} = \frac{1}{1 + \epsilon \tau} X$.

In our examples $X = \; {\rm grad} \; u$, and $u$ is constant on
hypersurfaces, so we have $X = \dot{u} \; N$.
We need a flow $\psi_{\tau}$ such that
\[
\frac{d}{d \tau} \left(f \circ \psi_{\tau}(m) \right) =
Y_{\tau}(\psi_{\tau}(m))f
\]
for all $f \in C^{\infty}(\overline{M})$.
For our choice of $Y_{\tau}$, the flow on ${\rm int}(I) \times P$
is of the form $\psi_{\tau}(t, p)=(\Xi(\tau, t), p)$ where
\[
\frac{d \Xi}{d \tau} = \frac{\dot{u}(\Xi)}{1 + \epsilon \tau}.
\]
Hence $\Xi$ is given by
\[
\Xi(\tau, t) = \left\{
\begin{array}{ll}
F^{-1} \left( \frac{\log(1 + \epsilon \tau)}{\epsilon} + F(t)\right)
    & \mbox{if $\epsilon \neq 0$}  \\
F^{-1} (\tau + F(t)) & \mbox{if $\epsilon =0$},
\end{array}
\right.    \]
where $F(t)$ is an antiderivative of $\frac{1}{\dot{u}(t)}=
\frac{1}{\kappa_1 f(t)}$. On the other hand, at an endpoint of
the interval $I$, $\dot{u}$ must be zero by Lemma \ref{gradRS}. So at
a compactifying submanifold of our shrinking or expanding solitons,
the Ricci flow just homothetically shrinks the submanifold as
we approach the critical time. This is consistent with $F(t)$
approaching $\pm \infty$ as $t\rightarrow 0,$ as we shall see below.

Let us consider our noncompact shrinking solitons. Now $\epsilon =-1$ and
$\kappa_1$ is positive, and $F(t)$ is asymptotically a positive constant
times $\log t$ as $t$ approaches $0$ or $\infty$. More precisely,
the discussion before Thm \ref{FIK} shows
that $F(t) \sim 2 \log t$ as $t$ becomes large, since $f(t) \sim
(2 \kappa_1)^{-1} t$. Hence $F^{-1}(t) \sim e^{\frac{1}{2}t}$ for $t$ large.
 As $\tau$ approaches $1$ from below, we see
\[
\Xi(\tau, t) \sim e^{\frac{1}{2}F(t)} (1 - \tau)^{-\frac{1}{2}}
\sim t (1-\tau)^{-\frac{1}{2}}.
\]
In particular $\Xi(\tau, t) \rightarrow \infty$ as $\tau$ tends to $1$.
Now the Ricci flow is the combination of $\psi_{\tau}$
and overall rescaling by $(1-\tau)$, so we see that the leading, i.e., $t^2$
terms in the metric coefficients survive and the other terms are killed.
So under the Ricci flow our soliton flows towards a cone, which is
the same as the asymptotic cone of the original metric.

For the complete expanding solitons, we have $\epsilon =1$ and
$\kappa_1 < 0$. Asymptotically
$\alpha(s) \sim -\frac{s}{\kappa_1}$ and $f(t) \sim (-2 \kappa_1)^{-1}t$,
so $F(t) \sim -2 \log t$ and tends to $-\infty$ for large $t$.
 Hence $F^{-1}(t) \sim e^{-\frac{1}{2}t}$
for $t$ large negative. We have
\[
\Xi(\tau, t)= F^{-1} \left(F(t) + \log(1 + \tau) \right),
\]
so as $\tau$ approaches $-1$ from above
\[
\Xi(\tau, t) \sim e^{-\frac{1}{2} F(t)} (1 + \tau)^{-\frac{1}{2}}
      \sim t (1+\tau)^{-\frac{1}{2}}.
\]
Under the Ricci flow this term is rescaled by $1 + \tau$, and again all
terms except the leading conical $t^2$ terms disappear in the limit.
So as before the soliton approaches the asymptotic cone of the original metric.

For the complete steady soliton, we have $\epsilon =0$ and $\kappa_1 < 0$. Now
$\alpha(s)$ and $f(t)$ are asymptotic to positive constants as $t$ becomes large,
while $f(t) \sim t$ for $t \sim 0$. Hence $F(t) \sim \mu t$ for $t$ large,
where $\mu$ is a negative constant; also $F(t) \sim \frac{\log t}{\kappa_1}$
for $t$ close to zero. Hence $F(t)$ tends to $\infty$ as $t \rightarrow 0$ and
tends to $-\infty$ as $t \rightarrow \infty$. So $F^{-1}(t) \sim \mu^{-1} t$
for $t \rightarrow -\infty$ and $\sim e^{\kappa_1 t}$ for $t \rightarrow
\infty$.

Hence as $\tau \rightarrow -\infty$, we have $\Xi(\tau, t) \sim \mu^{-1}
(\tau + F(t))$. In particular $\Xi(\tau, t) \rightarrow \infty$.

As $\tau \rightarrow + \infty$, we have $\Xi(\tau, t) \sim e^{\kappa_1 \tau}
\cdot e^{\kappa_1 F(t)}$, so $\Xi(\tau, t)$ tends to zero.

\begin{rmk}
Recall from the discussion after Eq.(\ref{mu-i}) that (\ref{mu-i})
also admits solutions where $\beta_i$ are quadratic polynomials in $s$.
In \cite{WW} such a choice of $\beta_i$ was shown to give rise to
Hermitian, non-K\"ahler Einstein metrics. However, for the Ricci soliton
equations such an ansatz turns out to be inconsistent except in the case
of a trivial soliton.
\end{rmk}

\section{\bf Another class of examples}

We may find another class of K\"ahler-Ricci solitons by exploiting the ideas
used to study Einstein metrics in \cite{DW}. Recall that in that paper we took
the hypersurface to be a homogeneous space $G/K$ (for $G$ compact and semisimple)
that fitted into a fibration
\begin{equation} \label{fibration}
S^1 = Q/K \rightarrow G/K \rightarrow G/Q
\end{equation}
where $G/Q$ is a generalised flag variety with a fixed invariant complex
structure (i.e., $G/Q$ is a connected compact homogeneous K\"ahlerian space
or coadjoint orbit for $G$).

In fact every circle bundle over $G/Q$ is of this form. For generic choices
of circle bundle the isotropy representation for $G/K$ is multiplicity free,
that is, all the irreducible summands are inequivalent as $K$-modules. We
shall always make this genericity assumption in the following discussion.
We write, therefore, the isotropy representation as:
\begin{equation} \label{isorep}
\p = \p_0 \oplus \p_1 \oplus \cdots \oplus \p_r
\end{equation}
where $\p_0$ is the 1-dimensional trivial represention corresponding to the
tangent space to $S^1$ in (\ref{fibration}). We denote by $d_i$ the real dimension
of $\p_i$. As $K$ acts trivially on $\p_0$, the adjoint action of $\p_0$ is
$K$-equivariant, and  it preserves each $\p_i, 1\leq i \leq r,$
since $\p$ has no multiplicities. We can choose a $G$-invariant complex structure
$J^*$ on $G/Q$, which on each $\p_i, \; 1 \leq i \leq r,$ is proportional
to ad$(X_0)$ for some $X_0 \in \p_0$. Note that each dimension $d_i$ is even.

The $G$-invariant metric $g_t$ on $G/K$ may now be written in the form
\begin{equation} \label{gt}
g_t = f(t)^2 \langle \,\, ,\,\, \rangle |_{\p_0} \oplus g_1(t)^2 \langle \,\, , \,\,\rangle |_{\p_1} \oplus
\cdots \oplus g_r(t)^2 \langle \,\, ,\,\, \rangle |_{\p_r}.
\end{equation}
Here the background metric $\langle \, , \, \rangle$ on $G/K$ is chosen to submerse
over the canonical K\"ahler-Einstein metric $g_{KE}^*$ on $G/Q$ with Einstein
constant equal to one. Moreover, if $U$ is the element of $\p_0$ such that
$\exp(2 \pi t U) = e^{2 \pi i t} \in S^1 = Q/K$, then we choose the background metric
on the $S^1$ fibre so that $\langle U, U \rangle = 1$. Observe from the above remarks
about the ad$(\p_0)$ action that $g_t$ is a Riemannian submersion over
the metric $g_{t}^* = \sum_{i=1}^{r} \, g_i(t)^2 \langle \,\, ,\,\, \rangle |_{\p_i}$
on $G/Q$. We denote the Ricci tensor of $g_{t}^{*}$ by Ric$_{t}^{*}$.

We denote by $\Omega^*$ the invariant 2-form on $G/Q$ defined by
\begin{equation} \label{Omega}
\Omega^*(Y,Z) = -\langle U, [Y,Z]|_{\p} \rangle.
\end{equation}
Now $-\frac{1}{2 \pi} \, \Omega^*$ represents the Euler class of the fibration
(\ref{fibration}). By invariance and the multiplicity free property,
$\Omega^*(\p_i, \p_j)=0$ for $i \neq j$.
Moreover there exist constants $b_i \; (i=1, \ldots, r)$ such that
\begin{equation} \label{OmegaKE}
\Omega^* |_{\p_i} = b_i \, \Theta_{KE}^* |_{\p_i} \;\;\; (i=1, \ldots, r),
\end{equation}
where $\Theta_{KE}^*$ is the K\"ahler form for $g_{KE}^*$.

As in \cite{DW} we can define an integrable complex structure $\bar{J}$ on
${\rm int}(I) \times (G/K)$ by lifting $J^*$ to the horizontal space in
(\ref{fibration}) and defining $\bar{J}(\frac{\partial}{\partial t})
= f^{-1} U$. Now $\bar{J}$ is Hermitian with respect to $\bar{g} = dt^2 + g_t$.

We now consider the Ricci soliton equations for the cohomogeneity one metric $\bar{g}=dt^2
+ g_t$ and the 1-form $\bar{\omega} = \dot{u} \; dt$,
where $u$ is $G$-invariant. By Prop 3.18 of \cite{BB},
$\overline{\ric} (X, N)=0$ due to the multiplicity free assumption, so as in
Remark \ref{commondecomp} there is no loss of generality in taking the soliton
to be of gradient type.

The equation (\ref{GRS-XN}) corresponding to
 mixed directions is now automatically satisfied,
so we just have to consider (\ref{GRS-NN}) and (\ref{GRS-XY}). We may
write the tensors in (\ref{GRS-XY}) as endomorphisms with respect to $g_t$, so the equation
becomes:
\begin{equation}
r_t - \dot{L} -({\tr L}) L + \dot{u} L + \frac{\epsilon}{2} =0
\end{equation}
where $r_t$ is the endomorphism defined by $\ric_t (X,Y) = g_t (r_t(X), Y)$.
Since $\p$ has no multiplicities, Schur's lemma implies that both $L$ and $r_t$
are diagonal with respect to (\ref{isorep}) and are scalar on each summand,
so (\ref{GRS-XY}) just becomes a system of $r+1$ scalar equations, one for each $\p_i$.

 We have
\[
L = {\rm diag} \left( \frac{\dot{f}}{f}, \frac{\dot{g_1}}{g_1} I_{d_1 \times d_1}, \cdots,
\frac{\dot{g_r}}{g_r} I_{d_r \times d_r} \right),
\]
where $I_{d_i \times d_i}$ denotes the $d_i \times d_i$ identity matrix.
Hence (\ref{GRS-NN}) is  just equation (\ref{NN}) with $2n_i$ replaced by $d_i$,
the dimension of $\p_i$. The component of (\ref{GRS-XY}) corresponding to $\p_0$ is
\begin{equation} \label{GRS-UU}
\frac{\ddot{f}}{f} + \sum_{i=1}^{r} d_i \frac{\dot{f} \dot{g_i}}{f g_i} -
\frac{\dot{u} \dot{f}}{f} - r_t = \frac{\epsilon}{2}
\end{equation}
where $r_t$ is the scalar defined by $\ric_t(U,U) = r_t g_t(U,U) = r_t f^2$.

Now, as in \cite{DW}, we may use the O'Neill formulae (\cite{Be} Chapter 9)
to compute $\ric_t(U,U)$. The submersion (\ref{fibration}) has totally geodesic fibres
so the O'Neill tensor $T$ is zero, and by (\cite{Be} 9.36a)
 $\ric_t(U,U) = g_t(AU, AU)$ where $A$ is the second O'Neill tensor, defined by
\[
A_{E_1} E_2 = {\mathcal H} \nabla^t_{{\mathcal H}E_1} {{\mathcal V}E_2} +
              {\mathcal V} \nabla^t_ {{\mathcal H}E_1} {{\mathcal H}E_2}.
\]
where $\mathcal H$ and $\mathcal V$ denote horizontal and vertical components respectively.

Let us take a basis $Y_{i \alpha} : \alpha = 1, \ldots, d_i$ for each $\p_i \;(i \geq 1)$,
orthonormal with respect to the background metric $\langle \, , \, \rangle$.
Now we form a $g_t$-orthonormal basis for the horizontal space $T (G/Q)$
by taking $Y_{i \alpha} / g_i : \alpha=1,\ldots, d_i, \;\; i =1, \ldots, r$.
We denote these vectors by $X_k$ and note that $A_{X_k} U$ is horizontal.
By \cite{Be} 9.36a and 9.33c, we have
$$\ric_t(U,U) = g_t (AU, AU) := \sum_{k} g_t (A_{X_k}U,A_{X_k}U). $$
It follows that
\begin{eqnarray*}
\ric_t(U, U) &=& \sum_{j,k} g_t(A_{X_k} U, X_j)^2 \\
            &=& \sum_{j,k} g_t(A_{X_k} X_j, U)^2 \\
            &=& \sum_{j,k} \frac{f^4}{4} \Omega^*(X_j,X_k)^{2}\\
            &=& \frac{f^4}{4}\sum_{i=1}^{r} \sum_{\alpha, \beta=1}^{d_i} \frac{b_i^2
 \Theta_{KE}^*(Y_{i \alpha}, Y_{i \beta})^2}{g_i^4} \\
            &=& \frac{f^4}{4} \sum_{i=1}^{r} \frac{d_i b_i^2}{g_i^4}
\end{eqnarray*}
where we have used \cite{Be} 9.21d, 9.24 and the definition (\ref{Omega})
and property (\ref{OmegaKE}) of $\Omega^*$.

Hence $r_t$ in equation (\ref{GRS-UU}) is $\sum_{i=1}^{r}
\frac{f^2}{g_{i}^{4}} \frac{d_i b_{i}^{2}}{4}$, and (\ref{GRS-UU}) is
equivalent to (\ref{UU}) with $d_i = 2n_i$ as above, and $b_i^2 =
q_i^2$.

The component of (\ref{GRS-XY}) corresponding to $\p_i \;(i \geq 1)$ is
\begin{equation} \label{GRS-ii}
\frac{\ddot{g_i}}{g_i} - \left( \frac{\dot{g_i}}{g_i} \right)^2 +
\frac{f \dot{g_i}}{fg_i}
+ \sum_{j=1}^{r} d_j \frac{\dot{g_i} \dot{g_j}}{g_i g_j} -
\frac{\dot{u} \dot{g_i}}{g_i}
- r_t = \frac{\epsilon}{2}
\end{equation}
where $r_t$ is the scalar defined by $\ric_t(Y_{i \alpha}, Y_{i \alpha})= r_t g_i^2$.

We take $\tilde{U} = \frac{1}{f} U$ as a $g_t$-orthonormal basis for the vertical space.
Now, using the O'Neill formulae again:
\begin{eqnarray*}
\ric_t(Y_{i \alpha}, Y_{i \alpha}) &=&  \ric_t^* (Y_{i \alpha}, Y_{i \alpha}) -
2 g_t(A_{Y_{i \alpha}} \tilde{U},  A_{Y_{i \alpha}} \tilde{U}) \\
                                  &=& \ric_t^*(Y_{i \alpha}, Y_{i \alpha}) -
2 \sum_{j} g_t(A_{Y_{i \alpha}} \tilde{U}, X_j)^2 \\
                                  &=& \ric_t^*(Y_{i \alpha}, Y_{i \alpha}) -
2 \sum_{j} g_t(A_{Y_{i \alpha}} X_j, \tilde{U})^2 \\
                                  &=&  \ric_t^*(Y_{i \alpha}, Y_{i \alpha}) -
2 \sum_{j} \frac{f^2}{4} \Omega^*(Y_{i \alpha},X_j)^2 \\
                                  &=& \ric_t^*(Y_{i \alpha}, Y_{i \alpha}) -
\frac{f^2}{2} \sum_{\beta =1}^{d_i} \frac{b_i^2 \Theta_{KE}^* (Y_{i \alpha}, Y_{i \beta})^2}{g_i^2} \\
                                  &=& \ric_t^*(Y_{i \alpha}, Y_{i \alpha}) -
\frac{b_i^2 f^2}{2 g_i^2}
\end{eqnarray*}
where we used \cite{Be} 9.36c, 9.33a as well as the earlier calculations.

Now we make the ansatz of \S3 (and of \cite{DW}); that is, we assume
that the functions $g_i^2$ are linear polynomials in $s$, the
antiderivative of $f(t)$, We further assume that $b_i$ is the
coefficient of $s$, so
\begin{equation}
g_i^2 (t) = b_i s + a_i  \;\;\;\  i=1, \ldots, r.
\end{equation}
We showed in \S 1 of \cite{DW} that this means that the
 2-form $\bar{\Theta}$ defined by the metric $\bar{g}$ and complex structure $\bar{J}$
 is closed in mixed directions,
i.e., $d \bar{\Theta} (\frac{\partial}{\partial t}, \cdot , \cdot )=0$. Also
$d \bar{\Theta} (U, \cdot, \cdot)=0$ so $\bar{\Theta}$ satisfies
the K\"ahler condition provided that
\begin{equation} \label{tau}
\Theta_{s}^{*} = \sum_{i=1}^{r} (b_i s + a_i) \Theta_{KE}^{*} |_{\p_i}
\end{equation}
is closed as a 2-form on $G/Q$ for all $s$. From (\ref{OmegaKE}) this is true
provided we choose $\sum_{i=1}^{r} a_i \Theta_{KE}^{*} |_{\p_i}$ to be closed.

If we choose $a_i$ in this way then for all $s$ the metric $g_s^*$ on $G/Q$ is
a $G$-invariant K\"ahler metric with respect to the fixed complex structure
$J^*$. As discussed in \cite{DW}, $g_{s}^*$ therefore has the same Ricci form
as $g_{KE}^*$, so
\[
\ric_t^*(Y_{i \alpha}, Y_{i \alpha}) = \ric_{KE}^* (Y_{i \alpha}, Y_{i \alpha})
= \langle Y_{i \alpha}, Y_{i \alpha} \rangle =1.
\]
Hence $r_t$ in (\ref{GRS-ii}) is
\[
\frac{1}{g_i^2} - \frac{b_i^2 f^2}{2 g_i^4}
\]
and now (\ref{GRS-ii}) is equivalent to  (\ref{XX}) with $p_i=1, b_i^2 = q_i^2$
and $d_i = 2n_i$.

\bigskip
So our equations are actually equivalent to those of \S 3, with the
ansatz, as in \S 3,  that each $g_i(t)$ is a linear function of $s$.
Putting $\alpha(s) = f(t)^2, \beta_i (s) = g_i(t)^2,
v = \prod_{i=1}^{r} {\beta_i}^{\frac{d_i}{2}}$ and $\phi(s) = u(t)$,
as in \S 3, the solution
is given by (\ref{alphaform1})
where
\begin{equation} \label{betaphi}
\beta_i(s) = b_i s + a_i  \;\;\ : \;\; \phi(s) = \kappa_1 (s + \kappa_0)
\end{equation}
and where we have consistency conditions
\begin{equation} \label{consist-3KE}
E^* = \frac{ \epsilon a_i + 2}{b_i}  \;\;\; : \;\;\; 1 \leq i \leq r.
\end{equation}
The constants $a_i, b_i$ are related to those in \S 3 by
$b_i = -q_i$ and $\sigma_i=\frac{a_i}{b_i}.$

If $\epsilon \neq 0$, then (\ref{consist-3KE}) and (\ref{tau}) show that
$\Theta_s$ is a linear combination of $\Omega$ and $\Theta_{KE}^*$, so
is automatically closed. If $\epsilon =0$, we do need to
impose the condition  that $\sum a_i \Theta_{KE}^* |_{\p_i}$ is closed.

The asymptotics of our solutions are the same as those in \S 3.
In the case of steady solitons with $\kappa_1 < 0$  the metric
is complete at infinity and the circle fibres in $g_t$ have asymptotically
constant radius. For expanding solitons with $\kappa_1 < 0$ again we have
completeness at infinity and the metric $\bar{g}$ is asymptotically conical.
We also have asymptotically conical shrinking solitons
with $\kappa_1 >0$.

As in \S 3, we can consider possible collapsing to special orbits.
Consider a special orbit $G/H$ where $K \subset H$ ; for smooth
collapsing we need $H/K$ to be a sphere $S^k$. As in \cite{DW}, we take
\[
\h = \kf \oplus \p_0 \oplus \p_1 \oplus \cdots \oplus \p_m
\]
for some $m \geq 0$ ($\h$ must be of this form
if the Hermitian structure extends over the special orbit).
 Note that $m=0$ corresponds to the case where
we just collapse a circle and we have $H=Q$. Also note that
$d_1 + \cdots + d_m = (k-1)$.

As before, for the metric to smoothly extend over the special orbit the
spheres in the normal bundle to $G/H$ must approach the round metric.
Let us write the round metric of constant curvature 1 on $H/K$ as
\[
c_0^2 \langle \,\, , \,\, \rangle |_{\p_0} + c_1^2 \langle \,\, , \,\, \rangle |_{\p_1} + \cdots +
c_m^2 \langle \,\, , \,\, \rangle |_{\p_m}.
\]
(In fact $c_0=1$ by our choice of $U$.) If the special orbit occurs at $s=t=0$, we need
\[
f(0) = 0 \;\;\; : \;\;\; g_i(0)=0 \; (i=1, \ldots, m)
\]
\[
\dot{f}(0) = 1 \;\;\; : \;\;\; \dot{g}_{i}(0) = c_i \; (i=1, \ldots, m)
\]
for the metric to extend smoothly. Similar calculations  to those in \S 3
show these conditions become
\[
\alpha(0)=0 \;\;\; : \;\;\; \beta_i (0)=0  \; (i=1, \ldots, m)
\]
\[
\alpha^{\prime}(0) =2 \;\;\ : \;\;\; \beta_i^{\prime}(0)= 2c_i^2 \;
(i=1, \ldots, m)
\]
so we need
\begin{equation} \label{ab}
a_i =0, \;\;\; b_i = 2c_i^2 \;\;\; (i=1, \ldots, m)
\end{equation}
from the conditions on $\beta_i$. The condition on $\alpha^{\prime}$
is equivalent to
\begin{equation} \label{Ek}
E^* = k+1.
\end{equation}
so by (\ref{consist-3KE})
\begin{equation} \label{ci}
c_i^2 = \frac{1}{k+1} \;\;\; : \;\; (i=1, \ldots, m)
\end{equation}

In the compact case, where the interval $I$ is $[0,s_*]$, we also need
conditions at $s=s_*$. These conditions (in terms of $\alpha, \beta$)
are the same as those at $s=0$ (for a different
set of indices $i$, of course) except that the signs of the derivatives
are changed. If $\tilde{k}$ denotes the dimension of the collapsing sphere at $s=s_*$,
and if the indices of the corresponding
collapsing summands are labelled $m+j, \ldots, r$,
then we have $c_i^2 = 1/(\tilde{k} +1) \; : \; (i=m+j, \ldots, r)$
and we need
\[
s_* = k + \tilde{k} + 2.
\]

As in \S 4 of \cite{DW}, we can arrange that the above conditions are satisfied
for suitable choices of $H$, and thus obtain complete examples of shrinking,
expanding or steady solitons, as well as examples of compact shrinking
solitons (as in \S 3, the integral (\ref{Futaki}) vanishes for some choice of
$\kappa_1$). In the case of compact shrinking solitons with no blowing-down
these examples are included in those found by different methods in \cite{PS2}.

\begin{rmk}
  To compare with the results of \cite{DW} for K\"ahler-Einstein
  metrics, observe that the Einstein constant $\Lambda$ equals
  $-\frac{\epsilon}{2}$, so we recover the consistency relation (2.13)
  of \cite{DW}. Our $E^*$ is twice the constant $C$ of \cite{DW}, and
  relations (\ref{ab}), (\ref{Ek}) give the relations (3.1),(3.2) of
  \cite{DW}.  Note also that the volume $v$ in the current paper
  differs from the ``$v$'' in \cite{DW} by a factor of $f^2$.
\end{rmk}

\bigskip
{\noindent \bf Acknowledgements.} We would like to thank Huai-Dong Cao and
Peng Lu for their comments on an earlier version of the paper.

\end{document}